\documentclass[12pt]{amsart}
\usepackage{epsfig,color}

\theoremstyle{definition}

\def\fnum{equation} 
\newtheorem{Thm}[\fnum]{Theorem}
\newtheorem{Cor}[\fnum]{Corollary}

\newtheorem{Lem}[\fnum]{Lemma}

\newtheorem{Rem}[\fnum]{Remark}
\newtheorem{Pro}[\fnum]{Proposition}

\numberwithin{equation}{section}

\newcommand{\nn}{{\bf{n}}}

\newcommand{\Tr}{{\text{Tr}}}

\newcommand{\Hess}{{\text {Hess}}}

\def\RR{{\mathbf R}}

\def\TT{{\mathbf T}}

\newcommand{\dv}{{\text {div}}}
\newcommand{\e}{{\text {e}}}

\newcommand{\cL}{{\mathcal{L}}}
\newcommand{\cLp}{{\mathcal{L}^+}}

\newcommand{\cP}{{\mathcal{P}}}

\newcommand{\eqr}[1]{(\ref{#1})}

\title{Level Set Method; partial regularity in $\RR^3$}

\author{Tobias Holck Colding}%
\address{MIT, Dept. of Math.\\
77 Massachusetts Avenue, Cambridge, MA 02139-4307.}
\author{William P. Minicozzi II}%

\thanks{The authors were partially supported by NSF DMS Grants 2405393 and 2304684.}

\email{colding@math.mit.edu and minicozz@math.mit.edu}

\begin{document}

\begin{abstract}
In the level set method, all level sets of a function move by mean curvature. The function is then a solution of a degenerate parabolic equation, the level set equation. Since the beginning of the 1990s it has been known that solutions are Lipschitz and, in general, not $C^2$. We show that, along a level set that is initially a smooth closed surface in $\RR^3$, the Hessian of the solution is bounded up to its first singular time. In fact, at each non-generic singularity the Hessian tends to zero, at a rate given by its instability. Using these bounds, we show in the sequel that, surprisingly, at the singular time the solution is twice differentiable and has a backward time derivative at every point of the level set, whatever the singularity, and identify both
the Hessian and time derivative.
\end{abstract}
\maketitle

\section{Introduction}

The level set method tracks a moving interface as the zero level set of a function, and moves every other level set of the function by the same law. It was introduced by Osher and Sethian for fronts that propagate with curvature dependent speed, \cite{OS}, cf. \cite{S}. Normals and curvatures are computed from the function, and changes of topology need no special treatment. This has made the method a standard tool in applied mathematics, from fluid dynamics and materials science to image processing and computer vision; see \cite{OF,GFO} and the references there. For motion by mean curvature, the function $u(x,t)$ solves
\begin{align}	\label{e:lse}
u_t=|\nabla u|\,\dv\,\frac{\nabla u}{|\nabla u|}\,  ,
\end{align}
with initial condition $u_0$, and each level set $\{u(\cdot,t)=c\}$ moves by mean curvature in a generalized sense. The analytic theory was developed by Evans and Spruck, \cite{ES1,ES2,ES3,ES4}, and independently, for a class of geometric equations, by Chen, Giga and Goto, \cite{ChGG}, \cite{GL}. It is set in the framework of viscosity solutions of Crandall and Lions, \cite{CrL,CrIL}, where uniqueness for second order equations rests on the comparison principle of Jensen, \cite{J}. Viscosity solutions of \eqr{e:lse} exist, are unique and are Lipschitz, \cite{ES1,ChGG}, but they are not in general $C^2$, \cite{I1}, cf. \cite{H2,KS,Se}. When the initial hypersurface is mean convex, the levels move monotonically, $u$ can be taken to be $a(x)-t$ with $a$ the arrival time, and \eqr{e:lse} becomes the degenerate elliptic equation $|\nabla a|\,\dv\big(\nabla a/|\nabla a|\big)=-1$. There the regularity is understood: $a$ is twice differentiable, \cite{CM3}, though the second derivatives might not be continuous;  \cite{CM5} determines when it is $C^2$; see also \cite{CM4}. For background on the connections between the geometry of the flow and the level set equation, see \cite{CM4,Cicm,CMf}.

The level set equation belongs squarely to the calculus of variations: by the coarea formula, $\int|\nabla u|\,dx$ is the total area of the level sets, and \eqr{e:lse} moves each of them by the gradient flow of its area. Its integrand $|p|$ is homogeneous of degree one, so the problem is degenerate. Hilbert, \cite{Hi}, set out two steps for the calculus of variations: existence, once the notion of solution is suitably extended (his 20th problem), and then the regularity of these solutions (his 19th).\footnote{In 1900 Hilbert presented 23 problems to the ICM; the 19th, 20th and 23rd concern the calculus of variations. The 20th asks whether every regular variational problem has a solution, ``provided also if need be that the notion of a solution shall be suitably extended''; the 19th asks for the regularity of solutions; and the 23rd asks for the further development of the calculus of variations; cf. \cite{All, Alm, Ca, DL}.} For the level set equation,
 the first step is the theory of viscosity solutions, and the second is the question we address here and in the sequel \cite{CMtwice}.

We study the Hessian of $u$ for general flows of surfaces in $\RR^3$, along the zero level set $\Sigma_t$ up to its first singular time $T$. Throughout, we assume:
\begin{enumerate}
\item[($\star$)] $u_0\in C^{\infty}(\RR^{3})$ is constant outside a compact set, $\Sigma=\{u_0=0\}$ is a closed embedded surface, $\nabla u_0\neq0$ on $\Sigma$, and $\Sigma_t$, $0\leq t<T$, is the classical mean curvature flow with $\Sigma_0=\Sigma$ on its maximal interval of existence.
\end{enumerate}
Every smooth closed embedded surface $\Sigma$ arises in this way: the signed distance $d$ to $\Sigma$ is smooth near $\Sigma$ \cite{KP}, so $u_0=\psi\circ d$ satisfies ($\star$) for any smooth nondecreasing $\psi$ with $\psi(s)=s$ near $0$ and $\psi$ constant for $|s|\geq\delta$, with $\delta>0$ small.

For $t<T$, the zero level set of $u(\cdot,t)$ is $\Sigma_t$, $u$ is smooth near $\Sigma_t$ and $\nabla u\neq0$ there (Lemma~\ref{l:Rsetup}). The function $v=|\nabla u|^{-1}$, restricted to $\Sigma_t$, is a positive solution of the Jacobi equation $(\partial_t-\Delta - |A|^2)v=0$ along the flow;  $v$ measures the spacing between nearby level sets. The Hessian of $u$ on $\Sigma_t$ splits into three parts. In an orthonormal frame $e_1,e_2,\nn$ with $\nn=\nabla u/|\nabla u|$,
\begin{align}	\label{e:Rsplit}
	\Hess_u(e_i,e_j)=-\frac{A(e_i,e_j)}{v}\, ,\qquad \Hess_u(e_i,\nn)=-\frac{e_i(v)}{v^2}\, ,\qquad \Hess_u(\nn,\nn)=-\frac{v_s}{v^3}\, ,
\end{align}
where $\partial_s=v\,\nn$ is the derivative across the levels. We call these the tangential, mixed and double normal parts, and write $m=|\nabla^Tv|/v^2$ for the norm of the mixed part.

{\it The degeneracy.} Written out, \eqr{e:lse} is $u_t=\big(\delta_{ij}-u_iu_j/|\nabla u|^2\big)\,u_{ij}$. The coefficient matrix is the projection onto the tangent plane of the level set, so the equation diffuses along each level and not at all across the levels, and it is undefined where $\nabla u=0$. If $\psi$ is increasing, then $\psi\circ u$ is again a solution, with the same levels, \cite[p.~639]{ES1}, and
\begin{align}
	\nabla(\psi\circ u)=\psi'\,\nabla u\, ,\qquad\Hess_{\psi\circ u}=\psi'\,\Hess_u+\psi''\,\nabla u\otimes\nabla u\, .\notag
\end{align}
So a bound for the Hessian must depend on the Hessian of $u_0$, through $\kappa=\sup\|\Hess_{u_0}\|$, with $\|\cdot\|$ the operator norm. The three parts are affected very differently. The tangential and mixed parts are only multiplied by $\psi'$; they are determined by a single level set and its Jacobi field $v$. The double normal part also changes by $\psi''|\nabla u|^2$, which can be prescribed freely. It is the variation $v_s$ of the Jacobi field across the levels, and it is not determined by the flow of the zero level, or by $v$ on it. A bound must be carried along the flow from the initial data, through equations whose coefficients degenerate at the singular points, where $\nabla u$ tends to zero. This is the main difficulty. 

{\it Tangent flows.} Tangent flows at time $T$ are smooth embedded shrinkers \cite{I2}, of multiplicity one  \cite{BK}, and each end of a noncompact one is asymptotic to a regular cone or to a round cylinder \cite{BK,Wa}. A tangent flow is unique if it does not depend on the sequence of rescalings. Uniqueness is known for cylinders \cite{CM1}, also in an effective form that is uniform in the center of the rescaling \cite{CMq}, for closed shrinkers \cite{Sch}, and for asymptotically conical shrinkers \cite{ChS}. This leaves the shrinkers that have a cylindrical end but are not cylinders, where uniqueness has recently been announced; we refer to it as (U). 

\begin{Thm}	\label{t:Rintro}
Assume ($\star$) and (U). Then
\begin{align}
	\sup_{t<T}\ \sup_{\Sigma_t}\ |\Hess_u|<\infty\, .\notag
\end{align}
\end{Thm}

The bound is not explicit.   The proof gives more at the singular points whose tangent flows are not round spheres or cylinders, which we call non-generic, or exotic. They are the unstable singularities \cite{CM7}, and there are infinitely many kinds of them, far from classified. A key point for the optimal regularity in the sequel is that at every exotic singular point the Hessian and the time derivative of $u$ tend to zero, at a rate. This is where regularity meets stability: the instability of the tangent flow makes $v$ outgrow the curvature, and this forces the decay. The instability is measured by the lowest eigenvalue $\mu_1$ of the stability operator of the tangent flow (Section~\ref{ss:Rtangent}, \cite{CM7}): $\mu_1\le-1$ for every shrinker other than a plane, with equality exactly for round spheres and cylinders, so $\mu_1<-1$ at every non-generic singular point \cite{CM7}.

\begin{Thm}	\label{t:Rdecay}
Under the hypotheses of Theorem~\ref{t:Rintro}, let $x_0$ be a non-generic singular point at time $T$, with $\mu_1$ as above, and let $\bar\rho=\max\{|x-x_0|,\sqrt{T-t}\}$. Then for every small $\theta>0$, on $\Sigma_t$ where $\bar\rho\leq r_1$, for some $r_1>0$,
\[
|\Hess_u|\leq C\bar\rho^{\,-2(\mu_1+1)-10\theta}\, ,\qquad |u_t|=\frac{|H|}{v}\leq C\bar\rho^{\,-2(\mu_1+1)-2\theta}\, ,
\]
with $C$ depending on $\theta$ and on the flow at a time $t_1 < T$.
\end{Thm}

The rate is twice the gap $-(\mu_1+1)>0$ between $\mu_1$ and $-1$, the eigenvalue of $H$: the more unstable the tangent flow, the faster the Hessian vanishes. At round spheres and cylinders the gap is zero, since exactly there $H$ is the lowest eigenfunction, cf.~\cite{CM7}, and the Hessian does not vanish \cite{CMtwice}. Theorem~\ref{t:Rdecay} follows from Theorems~\ref{new:mixed} and~\ref{new:dn}, with $2\alpha-1=-2(\mu_1+1)$ for the exponent $\alpha$ of Section~\ref{ss:Rcondhyp}.  
The mixed Hessian bound for $u$
   shows that $v$ is comparable to itself at its own scale.

\begin{Cor}	\label{c:Rharnack}
{\it Harnack inequality at the scale of $v$.} Under the hypotheses of Theorem~\ref{t:Rintro}, there is $C$ such that on $\Sigma_t$, for $t<T$,
\[
|\nabla^T\log v|\le C\,v\, ,\qquad |\partial_t\log v|\le C\,v^2\, ,
\]
where $\partial_t$ is the time derivative along the normal trajectories, as in the Jacobi equation. 
\end{Cor}

In a sequel \cite{CMtwice},
we show that, under the same hypotheses, $u(\cdot,T)$ is twice differentiable at every point of $\Sigma_T$, and identify its Hessian. This is surprising: unlike in the elliptic case, where all singularities are round spheres or cylinders, the singularities can be of infinitely many kinds, far from classified. It is also optimal, since $u(\cdot,T)$ is in general not $C^2$. Moreover $u$ has a backward time derivative. So $u$, a priori only a viscosity solution of \eqr{e:lse}, satisfies it classically at every point of $\Sigma_T$.

{\it The method.} All three parts of the Hessian are controlled by maximum principles along the flow $\Sigma_t$, and they are treated in order, each using only the earlier ones. The difficulty is at the singular points at time $T$, where $\nabla u=0$. All three parts are controlled by $v$ at the local geometric scale $\ell$. Here $\ell$ is a scale such that $|A|\le1/\ell$ on a parabolic ball of radius $\ell$; it is defined precisely in Section~\ref{ss:Rmixedpart}. Let $M$ be an upper bound for $v$ 
on this parabolic ball. Then
\[
\frac{|A|}{v}\le\frac1{\ell v}\, ,\qquad m\le C\,\frac Mv\cdot\frac1{\ell v}\, ,\qquad |\Hess_u(\nn,\nn)|\le C\Big(\frac Mv\Big)^2\frac1{\ell v}\, ,
\]
the second by interior estimates for the Jacobi equation (Lemma~\ref{new:gradcrit}), and the third since $\Hess_u(\nn,\nn)=-v_s/v^3$, once $|v_s|\le CM^2/\ell$. So it suffices to control two ratios, $1/(\ell v)$, which compares $v$ with the curvature, and $M/v$, which measures how well $v$ is known at the scale $\ell$, and to prove the bound for $v_s$. None of these is local: each comes from a maximum principle along the flow, started from the initial data.

The tangential part is the simplest: $|A|^2/v^2$ is a subsolution of a heat equation with a drift, so the maximum of $|A|/v$ on $\Sigma_t$ does not increase. In particular $v$ is at least of the order of the curvature.

The mixed part is $m=|\nabla^Tv|/v^2$, and the key is to control the Jacobi field $v$ from above and below. Section~\ref{s:Rtang} already bounds it from below by the curvature, $|H|\le2\kappa v$ (Lemma~\ref{l:Rnc}), and $v$ is bounded below by a positive constant.  Near generic singularities,
 the flow has mean convex neighborhoods, where $H>0$;  this lower bound is the right one, and what is missing is an upper bound by a constant times $H$. At  non-generic singularities, $v$ behaves like a Jacobi field of the self-similar flow, of order $(T-t)^{1/2+\mu_1}$ with $\mu_1 < -1$, while the curvature is of order $(T-t)^{-1/2}$. So $v$ outgrows the curvature, and we need a better lower bound and an upper bound that matches it up to a small power of the distance to the singular point.

We obtain both by comparing $v$ with barriers, built region by region:
\begin{itemize}
\item near the generic singularities, a multiple of $|H|$;
\item on the compact core of a non-generic singularity, the first eigenfunction of the stability operator, whose eigenvalue $\mu_1<-1$ measures the instability;
\item on its conical ends, powers of the distance;
\item on its cylindrical ends, which far out look like long necks, $|H|$ times such powers;
\item away from the singular set, a constant.
\end{itemize}
The barriers are built from the models. This is where uniqueness of the tangent flow is used: the flow has to stay close to a single model at all small scales. We glue these pieces into global super- and subsolutions of the Jacobi equation on $\Sigma_t$. The care is in the overlaps, where adjacent pieces must cross in the right order, and on the far parts of the cylindrical ends. The maximum principle then puts $v$ between a lower barrier and an upper barrier $B$, which differ only by small powers (Theorem~\ref{new:twosided}). Since $M\le CB$, this gives the bound for $m$ above, and near the non-generic points $1/(\ell v)\to0$ gives its decay.

The double normal part is $\Hess_u(\nn,\nn)=-v_s/v^3$, so the task is to bound $v_s$, the variation of $v$ across the levels, in terms of $v^3$. Differentiating the Jacobi equation across the levels shows that $v_s$ solves the same Jacobi equation, but with a source $S$ (Proposition~\ref{l:Rslave}). The source involves only the curvature, $v$ and its first two tangential derivatives, so the upper bound for $v$ and the interior estimates control it: $|S|\le CM^2/\ell^3$. We compare $v_s$ with a supersolution of size $B^2/\ell$ whose gain absorbs the source. Since the Jacobi operator does not involve $v$, it is built like $B$, each piece squared and divided by the local scale. Near the generic singularities, where $v$ is comparable to $|H|$, this is essentially $|H|^3$. The exception is near the tips of degenerate necks, where the flow is close to a translating bowl and $H^3$ fails to be a supersolution; there we use a lower power of $H$. Near the non-generic singularities,
 the pieces are larger powers of the same model barriers. Unlike $v$, $v_s$ is not determined by the flow of the zero level, so the comparison also needs a bound at the starting time. The maximum principle then gives $|v_s|\le CB^2/\ell$, and so the bound for the double normal part above.

Near the non-generic points, the two ratios have explicit sizes: $M/v\le C\bar\rho^{\,-4\theta}$ and $1/(\ell v)\le C\bar\rho^{\,-2(\mu_1+1)-2\theta}$, with $\bar\rho$ the parabolic distance of Theorem~\ref{t:Rdecay}. This gives the decay in Theorem~\ref{t:Rdecay}, and, with the analogous bound for the second derivatives of $v$, Corollary~\ref{c:Rharnack}.

Passages marked $\dagger$ concern only tangent flows that have a cylindrical end but are not cylinders. A reader willing to assume there are none can skip them.

\section{Preliminaries}	\label{s:Rprelim}

This section collects the facts about the level set flow, Jacobi fields, interior estimates and tangent flows that the three parts of the proof share, with proofs. The two used most are that $v$ is a Jacobi field (Lemma~\ref{l:Rsetup}) and the gradient criterion of Section~\ref{ss:Rtools}. 
 On $\Sigma_t$, $\nn=\nabla u/|\nabla u|$, $x_t=-H\nn$ with $H=\dv\,\nn$, $A(V,W)=\langle\nabla_VW,\nn\rangle$, so that $H=-\Tr_gA$, $v=|\nabla u|^{-1}$, $\partial_t$ is the derivative along the normal trajectories of the flow, $\partial_s=v\,\nn$ is the derivative across the levels, and $\kappa=\sup\|\Hess_{u_0}\|$, with $\|\cdot\|$ the operator norm. For a function $f>0$ we write $\cL_f=\Delta^T+\nabla^T_{\nabla^T\log f}$ and $\cLp=\cL_{v^2}$. We write $\mathcal S_T$ for the singular set at time $T$; a point of $\mathcal S_T$ is a neck or a sphere point if its tangent flow is a round shrinking cylinder or sphere.

We write $\Sigma_T=\{u(\cdot,T)=0\}$ for the level set at the first singular time. It is the set of limits of points of $\Sigma_t$ as $t\uparrow T$: such limits lie in it by continuity, and if $u\neq0$ on $B_r(x)\times(T-\delta,T)$, then $u$ has a fixed sign there, and comparison with a shrinking sphere keeps $|u(x,T)|>0$. 

\subsection{The zero level set}

\begin{Lem}	\label{l:Rsetup}
Assume ($\star$). For every $t<T$, the zero level set of $u(\cdot,t)$ is $\Sigma_t$, $u$ is smooth near $\Sigma_t$ and $\nabla u\ne 0$ on $\Sigma_t$. The restriction of $v$ to $\Sigma_t$ is a positive solution of the Jacobi equation \eqr{e:Jeqn} with $v=|\nabla u_0|^{-1}$ at $t=0$. For $T'<T$ the same holds, for $t\leq T'$, on the levels $\{u(\cdot,t)=c\}$ with $|c|\leq c_0(T')$.
\end{Lem}

\begin{proof}
The argument is standard; we include it. Fix $T'<T$. For $|c|$ small, the part of $\{u_0=c\}$ near $\Sigma$ is a normal graph over $\Sigma$, and these graphs and their classical flows on $[0,T']$ depend smoothly on $c$. Write them as $F(\cdot,t,c)$, with $F(\cdot,t,0)$ parametrizing $\Sigma_t$. Differentiating $u_0(F(\cdot,0,c))=c$ in $c$ gives $\langle F_c,\nn\rangle=|\nabla u_0|^{-1}>0$ at $t=0$. By Lemma~\ref{l:jfield}, $\langle F_c,\nn\rangle$ solves the Jacobi equation, so by \eqr{e:vminmono} it stays positive on $[0,T']$ for $|c|$ small. Hence $(y,t,c)\mapsto (F(y,t,c),t)$ is a local diffeomorphism, and it is injective, since distinct levels are disjoint by the avoidance principle and each level is embedded. The level set flow agrees with the classical flow as long as the latter is smooth, \cite{ES1,ChGG}, and the rest of $\{u_0=c\}$ stays away from it. Therefore $u(F(y,t,c),t)=c$ on the image, which is a neighborhood of $\bigcup_{t\leq T'}\Sigma_t\times\{t\}$. Differentiating in $c$ gives $\langle\nabla u,F_c\rangle=1$, that is, $\langle F_c,\nn\rangle=v$.
\end{proof}

\subsection{Jacobi fields and variations across the levels}	\label{ss:Rjacobi}

A Jacobi field is a solution of
\begin{align}	\label{e:Jeqn}
	\left( \partial_t - \Delta^T - |A|^2 \right) \, v = 0
\end{align}
along the   flow. If $\Sigma_t$ is closed and smooth for $t\in[t_0,t_1]$ and $v>0$ at time $t_0$, the maximum principle gives
\begin{align}	\label{e:vminmono}
	\min_{\Sigma_t}v\geq\min_{\Sigma_{t_0}}v>0\qquad\text{for }t\in[t_0,t_1]\, .
\end{align}
In particular $v\geq v_0$ on $\Sigma_t$ for $0\leq t<T$, where $v_0=1/\max_\Sigma|\nabla u_0|$.
For a tensor $\tau$ and $g\neq0$, the quotient rule $(\partial_t-\cL_{g^2})(\tau/g)=g^{-2}\big(g(\partial_t-\Delta^T)\tau-\tau(\partial_t-\Delta^T)g\big)$ holds, so the quotient of two Jacobi fields solves $(\partial_t-\cLp)(\tau/v)=0$ (see, e.g., Lemma~1.14 in \cite{CM1}). For a normal variation $F_s=v\,\nn$ of a family of immersions with metrics $g$, normals $\nn$ and second fundamental forms $A$, the variational formulas of \cite{HP} are
\begin{gather}
	g_s = - 2\, v \, A \, , \qquad \nn_s = - \nabla^T v \, , \qquad H_s = - (\Delta^T +|A|^2) \, v \, ,   \label{e:q1}  \\
	\partial_s |A|^2 = 2\, v \, \Tr_g A^3 + 2\, \langle A , \Hess^T_{v} \rangle \, .  \label{e:q5} 
\end{gather}

By \cite[Corollary 3.5]{H1}{\footnote{In \cite{H1} the second fundamental form is $h=-A$; the sign of $H$ is the same as here.}}, $H$ solves the Jacobi equation \eqr{e:Jeqn}, and
\begin{align}
	\left( \partial_t - \Delta^T \right) \, |A|^2 = 2 \, |A|^4 - 2\, |\nabla A|^2 \, .	\label{e:almostSimA2}
\end{align}

\begin{Lem}	\label{l:jfield}
If $F(x,t,s)$ is a mean curvature flow in $t$ for each fixed $s$, then $\langle F_s,\nn\rangle$ solves the Jacobi equation \eqr{e:Jeqn}.
\end{Lem}

\begin{proof}
By \eqr{e:q1} and the variation of $H$ under tangential motions, $H_s=-(\Delta^T+|A|^2)\langle F_s,\nn\rangle+\langle F_s^T,\nabla^TH\rangle$. Since $F_{st}=F_{ts}$, $F_t=-H\nn$, $\langle\nn_s,\nn\rangle=0$ and $\partial_t\nn=\nabla^TH$,
\begin{align}
	\partial_t\langle F_s,\nn\rangle=\langle(-H\nn)_s,\nn\rangle+\langle F_s,\nabla^TH\rangle=-H_s+\langle F_s^T,\nabla^TH\rangle=(\Delta^T+|A|^2)\langle F_s,\nn\rangle\, .\notag
\end{align}
\end{proof}

From now on $\partial_s=v\,\nn=\nabla u/|\nabla u|^2$ is the derivative across the levels. Differentiating $v=|\nabla u|^{-1}$ along it gives
\begin{align}	\label{e:vsformula}
	v_s=-\frac{\langle\nabla u,\nabla|\nabla u|\rangle}{|\nabla u|^4}=-v^3\,\Hess_u(\nn,\nn)\, .
\end{align}

\subsection{Interior estimates}	\label{ss:Rtools}

We write $\cP_r(x,t)=B_r(x)\times(t-r^2,t]$.

\begin{Lem}	\label{l:EH2}
If $\Sigma^n_t\subset\RR^{n+1}$ is a MCF of hypersurfaces and $v$ is a Jacobi field with $|A| , |v| \leq 1$ on $B_1 \times (0,1]$, then
\begin{align}	\label{e:EH2A}
	|\nabla^{k} A| + |(\nabla^T)^{k} v| \leq C_{k} {\text{ on }} B_{ \frac{1}{2}} \times [3/4,1] \, , 
\end{align}
where $C_{k}$ depends on $k$ and $n$.   
\end{Lem}

\begin{proof}
The bounds for $A$ are the interior estimates of   \cite[Theorem~3.4]{EH2}. Given them, there is $c=c(n)>0$ such that, near each point of $B_{1/2}\times[\frac34,1]$, each sheet of the flow is, on a parabolic cylinder of size $c$, a graph of bounded geometry over a fixed plane; several sheets may meet $B_1$, and they are treated separately. On each sheet $v$ solves $(\partial_t-\Delta^T
- |A|^2)v=0$, a linear parabolic equation whose coefficients are bounded with all their derivatives on these cylinders. So the bounds for $v$ follow from interior parabolic estimates.
\end{proof}

\begin{Lem}	\label{new:gradcrit}
{\it Gradient criterion.}
There is a universal constant $C$ with the following property. Let $(x,t)$ be a
point of the flow with $t<T$, and let $r>0$. Suppose that on $\cP_r(x,t)$ the flow, which may have several sheets, has $|A|\le1/r$, and that $v\le M$ there. Then at $(x,t)$
\[
|\nabla^Tv|\le\frac{CM}{r},\qquad |\Hess^T_v|\le\frac{CM}{r^2},\qquad
m\le\frac{CM}{r\,v^2}.
\]
The same holds on every level of $u$ near $\Sigma_t$.
\end{Lem}

\begin{proof}
Rescale by $1/r$ and apply Lemma~\ref{l:EH2}, with $n=2$, to $v/M$. This gives
the first two bounds, and the third follows since $m=|\nabla^Tv|/v^2$.
\end{proof}

{\it Barriers that are minima.} Several barriers below, in Lemma~\ref{l:Rsuper} and in Sections~\ref{s:Rmixed} and~\ref{s:Rdn},
are minima of finitely many smooth supersolutions. The maximum principle applies to such a minimum by a first-touching argument: at a first time at which a subsolution touches it from below, it touches the branch that realizes the minimum at that point, which is a smooth supersolution near it. In the same way, a maximum of finitely many smooth subsolutions is a
subsolution, as in Section~\ref{new:sec-lower}.

\subsection{Tangent flows}	\label{ss:Rtangent}

For a shrinker $S\subset\RR^3$, the stability operator  \cite{CM7} is
\[
L=\Delta-\tfrac12\langle y,\nabla\,\cdot\,\rangle+|A|^2+\tfrac12 {\text{ and }} L \, H = H \, .
\]

We will need the following properties of tangent flows:
\begin{Rem}	\label{p:Rstruct}
{\it Tangent flows and mean convex neighborhoods.}  Assume ($\star$).
\begin{enumerate}
\item Every tangent flow at a singular point at time $T$ is a smooth embedded shrinker of multiplicity one, and not a plane.
\item A singular point at time $T$ whose tangent flow is a round shrinking sphere or cylinder has a mean convex neighborhood.
\item Let $S$ be a tangent flow, $w$ a function on the self-similar flow $\sqrt{-t}\,S$, $t<0$, and $\hat w(y,s)=\sqrt{-t}\,w(\sqrt{-t}\,y,t)$, with $s=-\log(-t)$. Then $w$ is a Jacobi field if and only if $\partial_s\hat w=(L-1)\,\hat w$ on $S$.
\end{enumerate}
\end{Rem}

We explain these properties next.\\
(1) Let $(x_0,T)$ be a singular point, and let $\mathcal T$ be a tangent flow there: a Brakke limit of $\lambda_j(\mathcal M-(x_0,T))$ with $\lambda_j\to\infty$, where $\mathcal M$ is the flow $\Sigma_t$, $0\leq t<T$. Here $\mathcal M$ is a smooth flow of closed embedded surfaces, with bounded entropy and constant genus.

 By   \cite{I2} (see also \cite[Lemma 2.13(b)]{BK}), $\mathcal T$ is a smooth properly embedded shrinker with locally constant multiplicity. The multiplicity is constant, since shrinkers are connected \cite[Lemma 2.17]{BK}. The entropy bound gives polynomial volume growth \cite{CM1}.

  The flow is a bounded almost regular flow, and $\mathcal T$ is one of its blow-up limits, based at the final time. By the multiplicity one theorem   \cite[Theorems 1.1 and 1.2]{BK}, almost all tangent planes of $\mathcal T$ have multiplicity one. So the multiplicity is one.

 A multiplicity one plane has Gaussian density one, so $(x_0,T)$ would be regular by \cite{W1}.

(2) This was proven in Theorems $1.6$ and $1.7$ in  \cite{ChHH}. 

(3) Along normal trajectories of $\sqrt{-t}\,S$, a Jacobi field satisfies $\partial_tw=\Delta w+|A|^2w$. In the variables $(y,s)$, the normal trajectories move with velocity $-\tilde H\nn+\frac12y$, and on the shrinker $\tilde H=\frac12\langle y,\nn\rangle$. So the derivative at fixed $y$ picks up the transport term $-\frac12\langle y^T,\nabla\hat w\rangle$. Since $\partial_s$ of the factor $\sqrt{-t}$ contributes $-\frac12\hat w$, this turns the Jacobi equation for $w$ into 
 (3) for $\hat w$, and conversely; cf. \cite{CMcomplexity}.

\section{The tangential part}	\label{s:Rtang}

The tangential part $\Hess_u^{\TT}$ of the Hessian is $-A/v$, by \eqr{e:Rsplit}. It satisfies a maximum principle on its own, since $|A|^2/v^2$ solves a drift equation with a nonpositive right side (Lemma~\ref{l:evSi}). It also bounds $|A|$ and $|H|$ by $v$ (Lemma~\ref{l:Rnc}), which the other two parts use. Compare  Proposition~1.16 in \cite{CM1}; cf. subsection~$3.1$ in \cite{CIMW}.

 \begin{Lem}	\label{l:evSi}
 If $\Sigma_t$ is a mean curvature flow and $v> 0$ a Jacobi field, then
 \begin{align}	\label{e:evSi}
 	\left( \partial_t - \cLp \right) \, \frac{|A|^2}{v^2} & = - 2 \, \left| \nabla^T \frac{A}{v} \right|^2 \, .
 \end{align}
 \end{Lem}
 
 \begin{proof}
 Using the Jacobi equation  for $v$ and \eqr{e:almostSimA2} for $|A|^2$ gives
 \begin{align}
 	\left( \partial_t - \cLp \right) \, \frac{|A|^2}{v^2} &=  \frac{ \left( \partial_t - \cLp \right)   |A|^2}{v^2} - 
	\frac{|A|^2}{v^4} \, \left( \partial_t - \cLp \right)    v^2  + \frac{2}{v^4} \, \langle \nabla^T |A|^2 , \nabla^T v^2 \rangle - 2\, \frac{|A|^2}{v^6} \, |\nabla^T v^2|^2 \notag \\
	&= - 2 \, \frac{|\nabla A|^2}{v^2} + 2\, \frac{\langle \nabla^T v , \nabla^T |A|^2 \rangle}{v^3} - 2\, \frac{|A|^2}{v^4} \, |\nabla^T v|^2  = - 2 \, \left| \nabla^T \frac{A}{v} \right|^2 \, .  
 \end{align}
 \end{proof}

\begin{Pro}	\label{p:Rtang}
Assume ($\star$). For $t<T$,
\begin{align}
	\max_{\Sigma_t}\,\left|\Hess_u^{\TT}\right|\leq \max_{\Sigma}\,\left|\Hess_{u_0}^{\TT}\right|\, .
\end{align}
\end{Pro}

\begin{proof}
By Lemmas~\ref{l:Rsetup} and \ref{l:evSi}, $(\partial_t-\cLp)\,\frac{|A|^2}{v^2}\leq 0$ on $\Sigma_t$. Since $\Sigma_t$ is closed and smooth for $t\leq T'<T$, the maximum principle gives $\max_{\Sigma_t}|A|/v\leq\max_{\Sigma}|A|/v$, and $|A|/v=|\Hess_u^{\TT}|$ by \eqr{e:Rsplit}.
\end{proof}

By Proposition~\ref{p:Rtang}, $\Hess_u$ is bounded on $\Sigma_t$ up to $T$ if and only if the mixed part $|\nabla^T v|/v^2$ and the double normal part $|v_s|/v^3$ are. Recall that $m=\big|\nabla^T|\nabla u|\big|=|\nabla^Tv|/v^2$ is the norm of the mixed part, by \eqr{e:Rsplit}.

\begin{Lem}	\label{l:Rnc}
Assume ($\star$), and let $\kappa=\sup\|\Hess_{u_0}\|$. On $\Sigma_t$, for $t<T$, $|\Hess_u^{\TT}|\leq\sqrt2\,\kappa$, $|A|\leq\sqrt2\,\kappa\,v$ and $|H|\leq2\,\kappa\,v$.
\end{Lem}

\begin{proof}
On $\Sigma$, $\Hess_{u_0}^{\TT}$ is the restriction of $\Hess_{u_0}$ to a plane, of operator norm at most $\kappa$, so $|\Hess_{u_0}^{\TT}|\leq\sqrt2\,\kappa$. Proposition~\ref{p:Rtang} and $|A|/v=|\Hess_u^{\TT}|$ give the first two claims, and $|H|\leq\sqrt2\,|A|$ gives the third.
\end{proof}

\section{The mixed part}	\label{s:Rmixed}

Throughout this section and the next, $n=2$ and $(\star)$ holds. Here we bound $v$ above and below at the local scale and deduce the bound for the mixed part.

{\it The singular set.} Its generic part, the neck and sphere points, is compact away from small balls about the other singular points, so finitely many mean convex neighborhoods cover it (cf. Lemma \ref{new:finite}). On them $H>0$ controls the second fundamental form and its derivatives scale-invariantly, with uniform constants. The other singular points, whose tangent flows are not round spheres or cylinders, are finitely many (Lemma~\ref{new:finite}).

{\it The barriers.} The bounds come from the maximum principle for the Jacobi equation,  using
 barriers built from the geometry. The barriers are sub- and supersolutions modeled on
  Jacobi fields on the blowups and  powers of the distance on the ends. These pieces are glued together
  to  give global barriers on the closed surfaces $\Sigma_t$.  
 
 \vskip1mm
  The proofs of Lemmas~\ref{l:RNSnbhd}, \ref{l:Rss}--\ref{l:Rcompact} and~\ref{l:Rfarcyl} are in Appendix~\ref{app:geom}, and that of Lemma~\ref{new:polyb} in Appendix~\ref{new:appA}.

\subsection{Generic points}	\label{s:RNS}\label{ss:RNSnbhd}

Near  spherical and cylindrical singularities, 
 the flow has mean convex neighborhoods \cite{ChHH}. Lemma~\ref{l:RNSnbhd} records the structure of the flow on them: $H>0$, (MC) and canonical neighborhoods. The barrier there is $|H|$ itself, an exact Jacobi field, glued to $1$ away from the singular points (Section~\ref{new:sec-upper}). 

{\it Mean convex regions.} The functions $v$ and $H$ both solve the Jacobi equation $\partial_tw=\Delta w+|A|^2w$. For a set $\mathcal U$ of points $(x,t)$ with $x\in\Sigma_t$, such as the one in Lemma~\ref{l:RNSnbhd}, we say that (MC) holds on $\mathcal U$, with constants $h_0>0$ and $C_1$, if on $\mathcal U$
\begin{align}
	H\geq h_0\, ,\qquad |A|\leq C_1H\, ,\qquad |\nabla A|\leq C_1H^2\, ,\qquad |\nabla^2A|\leq C_1H^3\, .\notag
\end{align}
These  local estimates   follow from the canonical neighborhoods of Lemma~\ref{l:RNSnbhd}, cf. \cite{HK}.

\begin{Lem}	\label{l:RNSnbhd}
Let $U=\bigcup_{j=1}^{N_0}B_{r_j}(x_j)$, where the $x_j$ are generic singularities and each $B_{2r_j}(x_j)\times(T-4r_j^2,T+4r_j^2)$ is a mean convex neighborhood of $x_j$. Let $\delta\leq\min_jr_j^2$ and $\mathcal U=\{(x,t):x\in\Sigma_t\cap U,\ T-\delta\leq t<T\}$. Then:
\begin{enumerate}
\item $H>0$ on $\mathcal U$ for the orientation of the mean convex neighborhoods, which agree on overlaps since they are given by the direction of motion, and (MC) holds on $\mathcal U$, after decreasing $\delta$.
\item For every $\epsilon_0>0$ there is $\rho_0>0$ such that at every $(p,t)\in\mathcal U$ with $H(p,t)\geq1/\rho_0$, the flow rescaled by $H(p,t)$ about $(p,t)$ is $\epsilon_0$-close in $C^5$, on the backward parabolic ball of radius $1/\epsilon_0$, to a translating bowl or a round shrinking cylinder, or, near a sphere point, to a round shrinking sphere.
\end{enumerate}
\end{Lem}

Every point of $\mathcal S_T\cap\overline U$ is a neck or sphere point, since $\overline U\subset\bigcup_jB_{2r_j}(x_j)$, and the tangent flows at singular points in a mean convex neighborhood are mean convex shrinkers of multiplicity one other than planes, hence round spheres or cylinders \cite{CM7}.

{\it Orientation.} In Lemma~\ref{l:RNSnbhd}(1) the orientation is that of the mean convex neighborhoods, and it may differ from the one given by $\nabla u$.  A closed level set can have a neck of its inner region and a neck of its outer region. This does not matter. In the barriers $H$ means $|H|$, the mean curvature for the orientation of the mean convex neighborhoods, and the maximum principles of Sections~\ref{new:sec-upper}, \ref{new:sec-lower} and~\ref{new:sec-dn} are global on $\Sigma_t$ and bound $v$, which does not depend on the orientation, and $|v_s|$.

\subsection{Non-generic points}	\label{ss:Rcondhyp}\label{s:Rcond}\label{ss:Rlevel}\label{ss:Rcompact}

Throughout this subsection,  $x_0\in\Sigma_T$ is a non-generic singular point. The tangent flow is unique, by \cite{Sch}, \cite{ChS} or (U), and, by Remark~\ref{p:Rstruct}(1), the shrinker $S$ 
  has multiplicity one and is not a plane. By \cite{CM7}, the mean curvature
 $H_S$ changes sign,  
 the lowest eigenvalue $\mu_1$   of $L$ is  $<-1$, and  the lowest eigenfunction $\phi_1$ 
 with $L\phi_1=-\mu_1 \phi_1$ is positive. 
 Set
\begin{align}
	\alpha=-\mu_1-\tfrac12>\tfrac12\, ,\qquad y=\frac{x-x_0}{\sqrt\tau}\, ,\qquad s=-\log\tau\, ,\qquad \Phi(x,t)=\tau^{-\alpha}\,\phi_1(y)\, .\notag
\end{align}
Since $(L-1)\phi_1=(-\mu_1-1)\phi_1=(\alpha-\frac12)\phi_1$, the function $\sqrt\tau\,\Phi=\e^{(\alpha-\frac12)s}\phi_1(y)$ solves the equation in Remark~\ref{p:Rstruct}(3). So, by Remark~\ref{p:Rstruct}(3) after a translation in space and time, $\Phi$ is an exact positive Jacobi field on the self-similar flow $S_t=x_0+\sqrt\tau S$. 

\begin{Lem}	\label{new:polyb}
{\it Polynomial bounds.}
For each $\delta>0$ there are $c_\delta>0$ and $C_\delta,R_\delta$, depending only
on $S$ and $\delta$, such that on each end of $S$
\[
c_\delta|y|^{\gamma-\delta}\le\phi_1\le C_\delta|y|^{\gamma+\delta}
\qquad\text{for }|y|\ge R_\delta ,
\]
where $\gamma=-2\alpha$ on a conical end and $\gamma=1-2\alpha$ on a cylindrical
end.
\end{Lem}

Write $\rho=|x-x_0|$, and let $\tilde\Sigma_s=\tau^{-1/2}(\Sigma_t-x_0)$ be the rescaled flow.

\vskip1mm

\begin{Lem}	\label{l:Rss}
 Assume  $S$ is unique. There are $r_0>0$ and a nonincreasing function $\varepsilon(s)\to0$ such that, with $s_0=-\log r_0^2$ and $\hat R(s)=\big(\e^{-2/\varepsilon(s)}+\e^{-s/2}\big)^{-1/2}$:
\begin{enumerate}
\item[(i)] for $s\geq s_0$, $\tilde\Sigma_s\cap B_{\hat R(s)}$ is the normal graph over $S$ of a function $w(\cdot,s)$ whose $C^2$ norm on $S\cap B_{\hat R(s)}$ is at most $\varepsilon(s)$;
\item[(ii)] on a conical end with asymptotic cone $\mathcal C$, for $\hat R(s)\sqrt\tau\leq|x-x_0|\leq r_0$, the surface $\Sigma_t$ is the normal graph over $x_0+\mathcal C$ of a function $w$ with $|w|+\rho|\nabla w|+\rho^2|\nabla^2w|\leq\delta(\rho)\,\rho$, where $\delta(\rho)$ is nondecreasing and tends to $0$ as $\rho\to0$;
\item[(iii)] the static radius $\hat R(s)\e^{-s/2}$ is continuous and strictly decreasing, and $\varepsilon(s')\leq2\,\varepsilon(s)$ for $s-2\log\hat R(s)\leq s'\leq s$;
\item[(iv)] there are $R_G$, $c_g>0$ and $C_g$, depending only on $S$, such that at every point $(x,t)$ of a conical end with $R_G\sqrt\tau\le\rho\le r_0/2$, the flow on $\cP_{c_g\rho}(x,t)$ is a single sheet with $|\nabla^kA|\le C_g\rho^{-1-k}$ for $k\le2$;
\item[(v)] there is $\delta>0$ such that, for $T-\delta\le t<T$, $\Sigma_t\cap B_\delta(x_0)$ lies on the component of $\Sigma_t$ that contains the graph in (i).
\end{enumerate}
\end{Lem}

Here $\hat R$ is comparable to $\min\{\e^{1/\varepsilon(s)},\e^{s/4}\}$, and (iii) says that $\varepsilon$ varies slowly on the windows on which it is used. The entry time $s_\rho$ of the static radius $\rho$ is given by $\hat R(s_\rho)\e^{-s_\rho/2}=\rho$; it is continuous and decreasing in $\rho$, and we let $t_\rho=T-\e^{-s_\rho}$. In (ii), $\delta(\rho)$ depends only on $\varepsilon(s_\rho)$ and on how fast $S$ approaches its cone. On the far graphical part of a conical end, (i) alone gives only $|A|\le(C+C\varepsilon|y|)/\rho$, which is not scale invariant; (iv) gives the scale-invariant bound there, by a pseudolocality argument, cf.\ \cite{CIM,CM1}. We set $r_1=r_0/8$. Then the statements below that hold for $\rho\leq r_0/2$, and the definition of the far part $F$ of a cylindrical end, which uses the static radius $2\rho\leq r_0$, cover $\rho\leq4r_1$.

\begin{Lem}	\label{l:Rsheets}
{\it The flow near $x_0$.} Assume $S$ is unique. After decreasing $r_0$, for $\rho\leq r_0/2$ and $t_{2\rho}\leq t<T$, every point of $\Sigma_t$ with $\rho/2\leq|x-x_0|\leq\rho$ lies on the normal graph of Lemma~\ref{l:Rss}(ii) over a conical end, or within $C\rho/\hat R(s_{2\rho})$ of the axis of a cylindrical end.
\end{Lem}

{\it The core and the ends.} For $R_0\ge2$ we call $\{|y|\le R_0\}$ the core. The ends are the points with $|y|\ge R_0$ and $\rho\le r_0$; we take $R_0$ so large that $S\setminus B_{R_0}$ consists of the ends of $S$.

When the shrinker is compact, the component of $\Sigma_t$ near the point is close to it and there are no ends.

\begin{Lem} \label{l:Rcompact}
{\it Compact shrinkers.} If some tangent flow at $x_0\in\mathcal S_T$ is a compact shrinker $S$,
then   there exists $\delta>0$ such that, for $T-\delta\leq t<T$, Lemma~\ref{l:Rss}(i) holds, 
the graph there  is a  component $\mathcal C_t$ of $\Sigma_t$, and the other
components of $\Sigma_t$ do not meet $B_\delta(x_0)$.
\end{Lem}

\subsection{Barriers near a non-generic point}	\label{new:sec-radial}\label{ss:Rcondmixed}

 $\Phi$ is defined on $S_t$, and we transplant it to
$\Sigma_t$ by $\Phi(x,t)=\Phi(\pi_t(x),t)$, where $\pi_t$ is the nearest-point
projection onto $S_t$. We use this only on compact parts of the graphical region
below, where $\tilde\Sigma_s$ is a $C^2$-small normal graph over $S$ by Lemma~\ref{l:Rss}(i).

Each barrier is a positive Jacobi field of a model, times a factor that  absorbs the error from the flow not being the model. 
We use powers of the distance on the ends, with the factor $P_\pm$ (Lemma~\ref{new:radial}); on the core the model field is $\Phi$, with the factor $\tau^{\mp\gamma}$ (Lemma~\ref{new:core}); on the far part of a cylindrical end$^\dagger$ it is $H$ times a power of the distance, with a factor $P_\pm$ built from $1/H$, which absorbs the drift (Lemma~\ref{l:Rfarcyl}). The pieces are glued at the core by Lemma~\ref{new:glue}.

{\it Radial barriers.} On a conical end, the flow has bounded geometry at the scale of the distance $\rho$ to the singular point, and $v$ behaves like a power of $\rho$. Here we show that powers of $\rho$, times a factor $P_\pm$ that gains a little at the scale $\rho$, are super- and subsolutions of the Jacobi equation.

Fix $x_0\in\RR^3$, and write $\rho=|x-x_0|$, $\tau=T-t$ and
$|y|=\rho/\sqrt\tau$.

\begin{Lem}	\label{new:radial}
Let $R_0>0$ and $0<\eta\le\frac12$, and set
$P_\pm=1\pm\eta\bigl(1-R_0\sqrt\tau/\rho\bigr)$. On any mean curvature flow of
surfaces in $\RR^3$  at points with
$|y|\ge R_0$, and for $p>0$:
\begin{enumerate}
\item $(\partial_t-\Delta)\rho=(|\nabla^T\rho|^2-2)/\rho$;
\item if $\rho|A|\le C_g$ and $\eta R_0^2\ge4+6(p^2+C_g^2)$, then
$(\partial_t-\Delta-|A|^2)(P_+\rho^{-p})\ge\frac14\eta R_0|y|\,\rho^{-p-2}$;
\item if $\eta R_0^2\ge4+8p$, then
$(\partial_t-\Delta-|A|^2)(P_-\rho^{-p})\le-\frac14\eta R_0|y|\,\rho^{-p-2}$.
\end{enumerate}
Moreover $1\le P_+\le1+\eta$ and $1-\eta\le P_-\le1$ there.
\end{Lem}

\begin{proof}
Since $(\partial_t-\Delta)\rho^2=-4$, we have
$(\partial_t-\Delta)\rho=(|\nabla^T\rho|^2-2)/\rho$, which is (1). Hence
\[
(\partial_t-\Delta)\rho^{-p}=\rho^{-p-2}\bigl(2p-p(p+2)|\nabla^T\rho|^2\bigr),
\]
which lies between $-p^2\rho^{-p-2}$ and $2p\rho^{-p-2}$; for $p=1$ this is
$(\partial_t-\Delta)\rho^{-1}=\rho^{-3}(2-3|\nabla^T\rho|^2)$. Since
$\partial_t\sqrt\tau=-1/(2\sqrt\tau)$,
\[
\pm(\partial_t-\Delta)P_\pm=\eta R_0\Bigl(\frac{1}{2\sqrt\tau\,\rho}
-\frac{\sqrt\tau}{\rho^3}\bigl(2-3|\nabla^T\rho|^2\bigr)\Bigr)
\ge\eta\rho^{-2}\Bigl(\frac{R_0|y|}{2}-\frac{2R_0}{|y|}\Bigr),
\]
and $R_0/|y|\le1$, so $\pm(\partial_t-\Delta)P_\pm\ge\eta\rho^{-2}\bigl(\frac12R_0|y|-2\bigr)$. Moreover
\[
\nabla P_\pm=\pm\eta R_0\sqrt\tau\,\rho^{-2}\nabla\rho\, ,\qquad\nabla\rho^{-p}=-p\rho^{-p-1}\nabla\rho\, ,
\]
so
$\langle\nabla P_\pm,\nabla\rho^{-p}\rangle
=\mp p\eta R_0\sqrt\tau\rho^{-p-3}|\nabla^T\rho|^2$, and hence $\mp\langle\nabla P_\pm,\nabla\rho^{-p}\rangle\ge0$.
For (2), the product rule and these bounds give
\[
(\partial_t-\Delta-|A|^2)(P_+\rho^{-p})\ge\rho^{-p-2}
\Bigl[\eta\bigl(\tfrac12R_0|y|-2\bigr)-(1+\eta)(p^2+C_g^2)\Bigr],
\]
and $\frac12\eta R_0|y|\ge\frac14\eta R_0|y|+\frac14\eta R_0^2$, $2\eta\le1$,
$1+\eta\le\frac32$. For (3), drop $-|A|^2P_-\rho^{-p}\le0$ and
$-2\langle\nabla P_-,\nabla\rho^{-p}\rangle\le0$, and use
$P_-(\partial_t-\Delta)\rho^{-p}\le2p\rho^{-p-2}$. This gives the bound
$\rho^{-p-2}\bigl[-\eta(\frac12R_0|y|-2)+2p\bigr]$, and one argues as in (2).
\end{proof}

On the core $\{|y|\le R_0\}$, the flow is close to the shrinker, so $\Phi$ is nearly a Jacobi field. The error is small compared with $1/\tau$, so multiplying $\Phi$ by $\tau^{\mp\gamma}$ turns it into a super- or subsolution (Lemma~\ref{new:core}).
  On the graphical region, define
\begin{equation}	\label{e:Rpsieq}
E=-\frac{(\partial_t-\Delta-|A|^2)\Phi}{\Phi}.
\end{equation}
By
the interior estimates of \cite{EH2}, Lemma~\ref{l:Rss}(i) gives $C^k$ closeness on
$B_{\hat R(s)/2}$, with bounds $C_k\sup_{[s-1,s]}\varepsilon$. With Lemma~\ref{l:Rss}(iii), this gives for each $R\ge1$ a constant $C(R)$ such that, for
$s\ge s_0+1$ with $\hat R(s)\ge2R$,
\begin{equation}	\label{e:RE}
|E|\le C(R)\,\frac{\varepsilon(s)}{\tau}\qquad\text{on }\{|y|\le R\}.
\end{equation}

\begin{Lem}	\label{new:core}
{\it Barriers on the core.}
For $R\ge1$ and $\gamma>0$ there is $\varepsilon_C=\varepsilon_C(R,\gamma)>0$,
depending also on $S$, with the following property. If
$\varepsilon\le\varepsilon_C$ on $[s-1,s]$ and $\hat R(s)\ge2R$, then at time $s$,
on $\{|y|\le R\}$,
\[
(\partial_t-\Delta-|A|^2)\bigl(\tau^{-\gamma}\Phi\bigr)\ge
\frac{\gamma}{2\tau}\,\tau^{-\gamma}\Phi,\qquad
(\partial_t-\Delta-|A|^2)\bigl(\tau^{\gamma}\Phi\bigr)\le
-\frac{\gamma}{2\tau}\,\tau^{\gamma}\Phi .
\]
The same holds on the component $\mathcal C_t$ of Lemma~\ref{l:Rcompact} when $S$ is
compact, for $t$ close to $T$.
\end{Lem}

\begin{proof}
By \eqr{e:Rpsieq},
$(\partial_t-\Delta-|A|^2)(\tau^{\mp\gamma}\Phi)
=\tau^{\mp\gamma}\Phi\,(\pm\gamma/\tau-E)$, and $|E|\le C(R)\varepsilon/\tau$ by
\eqr{e:RE}. Take $\varepsilon_C=\gamma/(2C(R))$. 
\end{proof}

\begin{Lem}	\label{new:glue}
{\it Gluing at the core.} Let $B_c>0$ be a strict supersolution of the Jacobi equation on
$\{|y|\le2R_c\}$ and $B_e>0$ one on the ends $\{|y|\ge R_0\}$. Suppose that on each end
\[
\frac{B_e}{B_c}=F\,\frac{|y|^{a}}{\phi_1}\,,\qquad F_0^{-1}\le F\le F_0\,,\qquad a-\gamma<-\delta\,,
\]
with $\gamma$ the exponent of the end in Lemma~\ref{new:polyb}. There are $\kappa>0$ and $R_*$,
depending only on $F_0$, $R_0$, $\delta$, $S$ and the exponents $a$, such that if $R_c\ge R_*$,
the function $B$ equal to $\kappa B_c$ on $\{|y|\le R_0\}$, to $\min\{\kappa B_c,B_e\}$ on
$\{R_0\le|y|\le2R_c\}$ and to $B_e$ on $\{|y|\ge2R_c\}$ is continuous, and no solution of the
Jacobi equation touches it from below for a first time.

 For strict subsolutions with
$a-\gamma>\delta$ the same holds, with $\max$ in place of $\min$ and ``from above'' in place of
``from below''.
\end{Lem}

\begin{proof}
Let   $\omega_0$ and $\omega_1$ be the minimum and the maximum of $|y|^a/\phi_1$ on $\{R_0\le|y|\le2R_0\}$, over the ends. For supersolutions take $\kappa=\omega_0/(2F_0)$. Then $B_e\geq F_0^{-1}\omega_0B_c=2\kappa B_c$ on $\{R_0\le|y|\le2R_0\}$. Take $R_*\geq\max\{2R_0,R_\delta\}$. For $|y|\geq R_*$, Lemma~\ref{new:polyb} gives $B_e/(\kappa B_c)\le F_0(\kappa c_\delta)^{-1}|y|^{a-\gamma+\delta}<1$, once $R_*$ is large depending on $F_0$, $\kappa$ and $\delta$. So $B=\kappa B_c$ near $\{|y|=R_0\}$ and $B=B_e$ near $\{|y|=2R_c\}$, and $B$ is continuous. At a first touching point, the branch that realizes $B$ is a strict supersolution near it, which is impossible, as in the paragraph ``Barriers that are minima'' in Section~\ref{ss:Rtools}. For subsolutions take $\kappa=2F_0\omega_1$. Then $B_e\le\frac12\kappa B_c$ on $\{R_0\le|y|\le2R_0\}$, and $B_e/(\kappa B_c)\geq(F_0\kappa C_\delta)^{-1}|y|^{a-\gamma-\delta}>1$ for $|y|\geq R_*$, by the upper bound of Lemma~\ref{new:polyb}. The proof for subsolutions is similar.
\end{proof}

\medskip
{\it The far part of a cylindrical end.}$^\dagger$ The far part $F$ of a cylindrical end is the set of points of the flow within $C\max\{\sqrt\tau,\rho/\hat R(s_{2\rho})\}$ of the axis of the end with $R_0\sqrt\tau\leq\rho\leq4r_1$, where $C\geq2\sqrt2$ is the constant of Lemma~\ref{l:Rsheets}.

\begin{Lem}	\label{l:Rfarcyl}
{\it The far part of a cylindrical end.}$^\dagger$ Assume (U), with $S$ having a cylindrical end, and let $F$ be its far part. Let $\varepsilon_2>0$. If $r_1$ is small and $R_0$ is large, then:
\begin{enumerate}
\item The flow on $F$ is mean convex and, rescaled by $H$ about each point, $\varepsilon_2$-close in $C^5$, on the backward parabolic ball of radius $1/\varepsilon_2$, to a round shrinking cylinder or a translating bowl; so $|A|+|\nabla\log H|\leq C_AH$. Moreover $H\geq cR_0/\rho$ on $F$, with $c$ and $C_A$ depending only on $S$.
\item Let $0<\theta\le\frac14$, $P_\pm=1\pm\eta\bigl(1-cR_0/(\rho H)\bigr)$, with $c$ from (1), and $B_\pm=P_\pm\,\rho^{1-2\alpha\mp2\theta}$. There is $K_P$, depending only on $S$ and not on $\theta$, such that if   $K_P/R_0\le\eta\le\frac12$, then $1-\eta\le P_-\le1\le P_+\le1+\eta$ and, on $F$,
\[
\pm(\partial_t-\cL_{H^2})B_\pm\ge(1+4C_A)\,B_\pm\,\frac H\rho .
\]
So $HB_+$ is a supersolution and $HB_-$ a subsolution of the Jacobi equation on $F$.
\end{enumerate}
\end{Lem}

 The global barriers of the rest of this section and of Section~\ref{s:Rdn} involve constants that depend on each other. We fix the order in which they are chosen here, once for both sections.

Let $x^1,\dots,x^m$ be the non-generic singular  points at time $T$ (see Lemma \ref{new:finite}), with tangent flows $S^i$. At each $x^i$,
 as in
Section~\ref{s:Rcond},  $x_0=x^i$, $S=S^i$, $\alpha=\alpha_i$ and $\Phi=\Phi_i$, and we write
$\tau=T-t$, $\rho_i=|x-x^i|$, $y_i=(x-x^i)/\sqrt\tau$ and
$\bar\rho_i=\max\{\rho_i,\sqrt\tau\}$, which is the parabolic distance from $(x,t)$ to $(x^i,T)$ up to a factor $\sqrt2$. When $S^i$ is noncompact, $F^i$ is the (possibly empty)
union of the far parts of its cylindrical ends.

In the rest of this section, we bound $v$ above by a global barrier and below by local barriers near the $x^i$. In Section~\ref{s:Rdn} we bound $v_s$ above by a global barrier. Near the $x^i$ the bounds for $v$ are off
by factors $\bar\rho_i^{\,\pm2\theta}$, which is harmless by
Lemma~\ref{new:gradcrit}, since $\alpha_i>\frac12$ and $\theta\le(2\alpha_i-1)/12$ (see below).

Fix $\theta$ with
\[
0<\theta\le\theta_0:=\min\Bigl\{\frac14,\ \min_i\frac{2\alpha_i-1}{12}\Bigr\},
\]
and take $\delta=\theta$ in Lemma~\ref{new:polyb}. The other constants are fixed in order, each depending only on the $S^i$, $\theta$ and the constants before it; only the starting time $t_1$ and the overall factors of the barriers are fixed at the end, in Section~\ref{app:constants}. Fix $\varepsilon_2=\varepsilon'$, where $\varepsilon'>0$ is the constant of Lemma~\ref{l:Rsuper}, which depends only on the models. This fixes the constants $C_A$ and $c$ of Lemma~\ref{l:Rfarcyl}(1), and at each $x^i$ the constant $K_{P,i}$ of Lemma~\ref{l:Rfarcyl}(2), with $\alpha=\alpha_i$; let $K_P=\max_iK_{P,i}$, or $K_P=1$ if no $S^i$ has a cylindrical end. Next fix $R_0$ large, depending only on the $S^i$ and $\theta$: $R_0\ge\max\{2K_P,R_G,R_\delta\}$ (Lemmas~\ref{l:Rfarcyl}, \ref{l:Rss}(iv) and~\ref{new:polyb}), $S^i\setminus B_{R_0}$ consists of the ends of $S^i$, Lemma~\ref{l:Rfarcyl}(1) holds on $F^i$ with $\varepsilon_2$ once $r_1$ is small, and the finitely many further conditions of the form $R_0\ge C(S^i,\theta)$ stated where they are used in Sections~\ref{new:sec-upper}--\ref{new:sec-dn} hold. Let $\eta=K_P/R_0$, so that Lemma~\ref{l:Rfarcyl}(2) applies at each $x^i$. Then fix $r_0$ and $r_1=r_0/8$ so small that Lemma~\ref{l:Rss} holds at each $x^i$ on the pairwise disjoint balls $B_{r_0}(x^i)$, Lemma~\ref{l:Rsheets} holds for $\rho_i\le r_0/2$, Lemma~\ref{l:Rfarcyl}(1) holds with $\varepsilon_2$ on $F^i$, $r_0$ is below the $\delta$ of Lemma~\ref{l:Rcompact} at the $x^i$ with $S^i$ compact, and the conical ends of the $x^i$ lie within $c_g\rho_i$ of $x^i+\mathcal C$ where $|y_i|\ge R_0$ and $\rho_i\le4r_1$.

\subsection{The upper barrier for $v$}\label{new:sec-upper}

We build a supersolution $B_e$ away from the cores, glue it to the core barriers by Lemma~\ref{new:glue} into a supersolution $B$ on all of $\Sigma_t$, and apply the maximum principle (Proposition~\ref{new:vupper}).

{\it Radial profiles.} Fix $\zeta\in C^\infty(\RR)$ with $0\le\zeta\le1$, $\zeta=1$ on $(-\infty,r_1]$, $\zeta=0$ on $[2r_1,\infty)$, $|\zeta'|\le2/r_1$ and $|\zeta''|\le8/r_1^2$. For
$p>0$ let
\[
\psi_p(\rho)=\exp\Bigl(p\int_\rho^{2r_1}\zeta(\sigma)\,\frac{d\sigma}{\sigma}\Bigr).
\]
Then $\psi_p=a_p\rho^{-p}$ on $(0,r_1]$ for a constant $a_p>0$, $\psi_p=1$ on
$[2r_1,\infty)$, $|\rho(\log\psi_p)'|\le p$ and $\rho^2|\psi_p''|\le C(p)\psi_p$; by
Lemma~\ref{new:radial}(1) also $\rho^2|(\partial_t-\Delta)\psi_p|\le C(p)\psi_p$.

{\it The barrier $\Upsilon$.} Let $U_1=B_R\setminus\bigcup_i\bar B_{r_1}(x^i)$, with $\Sigma_t\subset B_{R/2}$. The singular points in $\overline{U_1}$ are neck and sphere points, and they form a compact set. Cover it by finitely many balls $B_{r_j/2}(x_j)$, with $B_{r_j}(x_j)$ as in Lemma~\ref{l:RNSnbhd} and so small that $B_{2r_j}(x_j)$ misses the conical ends of the $x^i$; this is possible since the axes of the cylindrical ends are at distance at least $c\rho_i$ from them. Let $U\supset U'$ be the unions of the $B_{r_j}(x_j)$ and of the $B_{r_j/2}(x_j)$. For $T-\delta_\Upsilon\le t<T$, (MC) holds on $\Sigma_t\cap U$ by Lemma~\ref{l:RNSnbhd}(1). After decreasing $\delta_\Upsilon$, the curvature and its first two derivatives are also bounded on $(\Sigma_t\cap U_1)\setminus U'$, since the points of $\Sigma_T\cap\overline{U_1}$ outside $U'$ are regular and points of $\Sigma_t$ accumulate as $t\to T$ only on $\Sigma_T$.

Let $V$ be the set of points within $c_gr_1$ of the axes of the cylindrical ends of the $x^i$ with $\frac32r_1<\rho_i<\frac72r_1$. Then $\Sigma_t\cap V\subset F^i$ for $t_1\le t<T$ (see below), so $H\ge cR_0/(4r_1)$ there by Lemma~\ref{l:Rfarcyl}(1). Fix $\chi\in C^\infty_c(U\cup V)$ with $0\le\chi\le1$, and $\chi=1$ on $U'$ and on $V\cap\{2r_1\le\rho_i\le3r_1\}$. Then $\chi=0$ on the conical ends, and on $\Sigma_t$ it varies only on $(U\cup V)\setminus U'$, where $H$ has a positive lower bound $h$ and the geometry is bounded. Let
\[
\Upsilon=\begin{cases}
H & \text{on } F^i\cap\{\rho_i\le3r_1\},\\
1 & \text{on the conical ends of the } x^i,\\
\chi|H|+1-\chi & \text{elsewhere.}
\end{cases}
\]
These agree where the cases meet, since $\chi=1$ on $F^i\cap\{2r_1\le\rho_i\le3r_1\}$ and $\chi=0$ on the conical ends. Where $\Upsilon=|H|$, it solves the Jacobi equation. Elsewhere on $U_1$, since $(\partial_t-\Delta^T)\chi=-\Tr_{\Sigma_t}D^2\chi$,
\[
(\partial_t-\Delta-|A|^2)\,\Upsilon=-(|H|-1)\Tr_{\Sigma_t}D^2\chi-2\langle\nabla\chi,\nabla|H|\rangle-(1-\chi)|A|^2
\]
is bounded, while $\Upsilon\ge\min\{1,h\}>0$. So $|(\partial_t-\Delta-|A|^2)\Upsilon|\le C_\Upsilon$ on $U_1$. 

{\it The set $V$.}$^\dagger$ We show that $\Sigma_t\cap V\subset F^i$ for $t_1\le t<T$. By the choice of $r_1$, the conical ends of $x^i$ lie within $c_g\rho_i$ of $x^i+\mathcal C$. This cone is at distance at least $2c_g\rho_i$ from the axes, since the sets $\Gamma_E$ in the proof of Corollary~\ref{c:BK1} are pairwise disjoint (after decreasing $c_g$, which keeps Lemma~\ref{l:Rss}(iv)). So the conical ends avoid $V$, and by Lemma~\ref{l:Rss}(i) before $t_{2\rho_i}$ and Lemma~\ref{l:Rsheets} after it, the other points of $\Sigma_t$ in $V$ lie near the axes, on $F^i$.

Away from the cores $\{|y_i|<R_0\}$ let
\begin{equation}\label{new:BU}
B_e=e^{C_4(t-t_1)}\,\hat P\,\hat\Psi\,\Upsilon ,
\end{equation}
where, near each $x^i$ with $S^i$ noncompact,
\[
\begin{array}{l|c|c}
 & \hat P & \hat\Psi\\ \hline
\text{conical ends of } x^i & 1+\eta-\eta R_0\sqrt\tau\,\zeta(\rho_i)/\rho_i & \psi_{2\alpha_i+2\theta}(\rho_i)\\
F^i & 1+\eta-\eta cR_0\,\zeta(\rho_i)/(\rho_iH) & \psi_{2\alpha_i-1+2\theta}(\rho_i)
\end{array}
\]
with $c$ from Lemma~\ref{l:Rfarcyl}(1), and $\hat P=1+\eta$, $\hat\Psi=1$ elsewhere. For $\rho_i\le r_1$, $\hat P$ is the $P_+$ of Lemma~\ref{new:radial} on the conical ends and of Lemma~\ref{l:Rfarcyl}(2) on $F^i$, so there $B_e=e^{C_4(t-t_1)}\,a\,P_+\rho_i^{-2\alpha_i-2\theta}$ and $B_e=e^{C_4(t-t_1)}\,a'HB_+$, with $a=a_{2\alpha_i+2\theta}$ and $a'=a_{2\alpha_i-1+2\theta}$. For $\rho_i\ge2r_1$, $B_e=e^{C_4(t-t_1)}(1+\eta)\Upsilon$. Since $\rho_iH\ge cR_0$ on $F^i$, $1\le\hat P\le1+\eta$ where $|y_i|\ge R_0$.

\begin{Lem}	\label{new:upperlemma}
If $C_4$ is large, depending on $r_1$ and on the flow away from the $x^i$, then $B_e$
is a strict supersolution of the Jacobi equation on
$\Sigma_t\cap\{|y_i|\ge R_0\ \text{for all}\ i\}$, for $t_1\le t<T$.
\end{Lem}

\begin{proof}
The factor $e^{C_4(t-t_1)}$ adds $C_4B_e$ to $(\partial_t-\Delta-|A|^2)B_e$. So it
suffices to show that the rest is at least $-C(r_1)B_e$, and positive where
$\rho_i\le r_1$.

{\it Conical ends, $\rho_i\le r_1$.} This is Lemma~\ref{new:radial}(2) with $p=2\alpha_i+2\theta$. Its hypothesis $\rho|A|\le C_g$ holds by
Lemma~\ref{l:Rss}(iv), since $|y_i|\ge R_0\ge R_G$, and $\eta R_0^2=K_PR_0\ge4+6(p^2+C_g^2)$ since $R_0$ is large.

{\it The far parts, $\rho_i\le r_1$.}$^\dagger$ On $F^i\cap\{\rho_i\le r_1\}$ we have $B_e=e^{C_4(t-t_1)}a'HB_+$, with $B_+$ from Lemma~\ref{l:Rfarcyl}(2), so this follows from that lemma.

{\it The transition, $r_1\le\rho_i\le2r_1$.} Let $Q=\hat P\hat\Psi$, and let $p$ be the exponent of $\hat\Psi$.

On the conical ends, $Q$ is a function of $\rho_i$ and $\tau$, with $-\partial_\tau\hat P\ge0$ and $|\partial_\rho\hat P|\le3\eta R_0\sqrt\tau/r_1^2$. So, by the choice of $\tau_1$,
\[
|\partial_\rho Q|\le\bigl(p+6K_P\sqrt\tau/r_1\bigr)\frac{Q}{\rho_i}\le(p+1)\frac{Q}{\rho_i},
\qquad\rho_i^2|\partial_\rho^2Q|\le C(r_1)\,Q .
\]
There $\Upsilon=1$, and the geometry is bounded at the scale $r_1$ by Lemma~\ref{l:Rss}(iv), so $(\partial_t-\Delta-|A|^2)Q\ge-C(r_1)Q$.

On $F^i$, $\hat P=1+\eta-\eta cR_0\,g/H$ with $g=\zeta(\rho_i)/\rho_i$. The computation of Step~4 of the proof of Lemma~\ref{l:Rfarcyl} applies with $g$ in place of $1/\rho$. Since $H\ge cR_0/(2r_1)$ there, the derivatives of $\zeta$ add only terms bounded by $C(r_1)Q$. Since $\eta R_0=K_P$, the gain $\eta cR_0\,\zeta|A|^2Q/(\rho_iH\hat P)\ge\eta cR_0\,\zeta HQ/(3\rho_i)$ absorbs the drift of $\hat\Psi$, which is at most $2pC_A\,\zeta HQ/\rho_i$ since $\partial_\rho\log\hat\Psi=-p\zeta/\rho_i$; both carry the factor $\zeta$. All other terms, including the cross term $-2\langle\nabla\hat P,\nabla\hat\Psi\rangle$, are at most $C(r_1)Q$, since $\rho_i\ge r_1$, $H\ge cR_0/(2r_1)$ and $|\nabla H^{-1}|\le C_A$.
So $(\partial_t-\cL_{H^2})Q\ge-C(r_1)Q$, and
$(\partial_t-\Delta-|A|^2)(QH)=H(\partial_t-\cL_{H^2})Q\ge-C(r_1)QH$.

{\it Elsewhere.} There $B_e=e^{C_4(t-t_1)}(1+\eta)\Upsilon$, and
$|(\partial_t-\Delta-|A|^2)\Upsilon|\le C_\Upsilon$.

Taking $C_4$ larger than the constants $C(r_1)$ and $C_\Upsilon$ proves the lemma.
\end{proof}

We will now define $B$ by blending $B_e$ with the separation of variables solutions on the core.
At each $x^i$ with $S^i$ noncompact, $\tau^{-\theta}\Phi_i$ is a strict supersolution on $\{|y_i|\le2R^i_c\}$, by Lemma~\ref{new:core} with $\gamma=\theta$. At fixed $|y_i|$, $B_e$ and $\tau^{-\theta}\Phi_i$ have the same power $\tau^{-\alpha_i-\theta}$. On a conical end
\[
\frac{B_e}{\tau^{-\theta}\Phi_i}=a\,e^{C_4(t-t_1)}P_+\,
\frac{|y_i|^{-2\alpha_i-2\theta}}{\phi_1},
\]
and on the graphical part of a cylindrical end
\[
\frac{B_e}{\tau^{-\theta}\Phi_i}=a'e^{C_4(t-t_1)}P_+\bigl(\sqrt\tau H\bigr)
\frac{|y_i|^{1-2\alpha_i-2\theta}}{\phi_1},
\]
where $1\le e^{C_4(t-t_1)}\le e$, $1\le P_+\le\frac32$ and $\frac12\le\sqrt\tau H\le1$. So Lemma~\ref{new:glue} applies, with $B_c=\tau^{-\theta}\Phi_i$, the $B_e$ above, $a-\gamma=-2\theta<-\delta$, and $F_0$ depending on $a$ and $a'$, hence on $r_1$. Let $K^i$ be its constant $\kappa$, and $B^i_c=K^i\tau^{-\theta}\Phi_i$. We fix one radius $R^i_c$, at least the $R_*$ of this gluing and of the gluings in Sections~\ref{new:sec-lower} and~\ref{new:sec-dn}, whose $R_*$ depend only on constants already fixed. Let $B$ be the function it gives: $B^i_c$ on $\{|y_i|\le R_0\}$, $\min\{B^i_c,B_e\}$ on $\{R_0\le|y_i|\le2R^i_c\}$, and $B_e$ elsewhere. At an $x^i$ with $S^i$ compact, let $B=B^i_c=\tau^{-\theta}\Phi_i$ on the component $\mathcal C_t$ of Lemma~\ref{l:Rcompact}; it is a strict supersolution by Lemma~\ref{new:core}.

\begin{Pro}	\label{new:vupper}
Let $N=\max_{\Sigma_{t_1}}(v/B)$. Then $v\le NB$ on $\Sigma_t$ for $t_1\le t<T$.
\end{Pro}

\begin{proof}
By Lemmas~\ref{new:upperlemma}, \ref{new:core}
and~\ref{new:glue}, $B$ is a positive continuous function that no solution of the
Jacobi equation touches from below for a first time. Since $\Sigma_t$ is closed for
$t<T$, the maximum principle on $[t_1,t']$, for each $t'<T$, gives $v\le NB$.
\end{proof}

\subsection{The lower barriers}\label{new:sec-lower}
Lower barriers are needed only near the non-generic points: by \eqr{e:vminmono} and Lemma~\ref{l:Rnc}, $v\geq\max\{v_0,|H|/(2\kappa)\}$ everywhere, which suffices elsewhere.

Fix $i$ with $S^i$ noncompact, and let
$D_i=\{(x,t):x\in\Sigma_t\cap B_{r_1}(x^i),\ t_1\le t<T\}$. Let
\begin{align*}
w_c&=\tau^{\theta}\Phi_i &&\text{on } \{|y_i|\le2R^i_c\},\\
w_e&=P_-\rho_i^{-2\alpha_i+2\theta} &&\text{on the conical ends},\\
w_f&=HP_-\rho_i^{1-2\alpha_i+2\theta} &&\text{on } F^i .
\end{align*}
These are strict subsolutions of the Jacobi equation: $w_c$ by
Lemma~\ref{new:core} with $\gamma=\theta$; $w_e$ by Lemma~\ref{new:radial}(3) with
$p=2\alpha_i-2\theta$, since $\eta R_0^2=K_PR_0\ge4+8p$ as $R_0$ is large; and
$w_f=HB_-$, with $B_-$ from Lemma~\ref{l:Rfarcyl}(2). On both kinds of ends the ratio of $w_e$ or $w_f$ to $w_c$ has the form of Lemma~\ref{new:glue}, with $F=P_-$ or $F=P_-\sqrt\tau H$ and $a-\gamma=2\theta>\delta$. Write $w_*$ for $w_e$ on the conical ends and for $w_f$ on $F^i$, and define $W$  from Lemma~\ref{new:glue}
\[
W=\begin{cases}
\kappa w_c & \text{on } \{|y_i|\le R_0\},\\
\max\{\kappa w_c,w_*\} & \text{on } \{R_0\le|y_i|\le2R^i_c\},\\
w_* & \text{on the rest of } D_i .
\end{cases}
\]
It is continuous, and no solution of the Jacobi equation touches it from above for a first time. At an $x^i$ with $S^i$ compact, let $W=\tau^\theta\Phi_i$ on $\mathcal C_t$.

\begin{Pro}	\label{new:vlower}
There is $c>0$ such that $v\ge cW$ on $D_i$, and on $\mathcal C_t$ for $t_1\le t<T$ when $S^i$ is compact.
\end{Pro}

\begin{proof}
The parabolic boundary of $D_i$ consists of the slice $t=t_1$ and of
$\Sigma_t\cap\partial B_{r_1}(x^i)$, which lies on the conical ends and on $F^i$,
since $2R^i_c\sqrt{\tau_1}\le r_1$. There $W=w_e\le r_1^{-2\alpha_i+2\theta}$ or
$W=w_f\le r_1^{1-2\alpha_i+2\theta}H$, while $v\ge v_0$ by \eqr{e:vminmono} and $v\ge H/(2\kappa)$ by Lemma~\ref{l:Rnc}. At $t=t_1$, $W\le C(t_1)$, since the flow is smooth there, and $v\ge v_0$. So $cW\le v$ on the parabolic boundary for $c$ small,
and the maximum principle gives $cW\le v$ on $D_i$. When $S^i$ is compact,
$\tau^\theta\Phi_i$ is a strict subsolution on $\mathcal C_t$ by
Lemma~\ref{new:core}, and $\mathcal C_t$ has no boundary.
\end{proof}

\subsection{The mixed part}	\label{ss:Rmixedpart}

 For ease of reference, we record the two-sided bounds on $v$:
 
\begin{Thm}	\label{new:twosided}
{\it Two-sided bounds.} There are $c>0$ and $C$ such that, for $t_1\le t<T$:
\begin{enumerate}
\item[(a)] at each $x^i$ with $S^i$ noncompact,
$c\,\tau^{\theta}\Phi_i\le v\le C\tau^{-\theta}\Phi_i$ on $\{|y_i|\le R^i_c\}$;
$c\,\rho_i^{-2\alpha_i+2\theta}\le v\le C\rho_i^{-2\alpha_i-2\theta}$ on the conical
ends with $R_0\sqrt\tau\le\rho_i\le r_1$; and
$c\,\rho_i^{1-2\alpha_i+2\theta}\le v/H\le C\rho_i^{1-2\alpha_i-2\theta}$ on
$F^i\cap\{\rho_i\le r_1\}$;
\item[(b)] at each $x^i$ with $S^i$ compact,
$c\,\tau^{\theta}\Phi_i\le v\le C\tau^{-\theta}\Phi_i$ on $\mathcal C_t$;
\item[(c)] elsewhere, $\max\{v_0,|H|/(2\kappa)\}\le v\le C\Upsilon$.
\end{enumerate}
\end{Thm}

\begin{proof}
The upper bounds follow from Proposition~\ref{new:vupper}. The lower bounds in (a) and (b) follow from
Proposition~\ref{new:vlower}. The lower bound in (c)  uses  \eqr{e:vminmono}
and   Lemma~\ref{l:Rnc}.
\end{proof}

We next define the local geometric scale $\ell$: For $t_1\le t<T$,  let  $\ell=\sqrt\tau$ on the cores $\{|y_i|\le2R^i_c\}$ and on the compact components
$\mathcal C_t$; $\ell=\rho_i$ on the conical ends with $\rho_i\le r_1$; $\ell=1/H$ on
$F^i\cap\{\rho_i\le r_1\}$ and where $\chi=1$; and $\ell=1$ elsewhere. 

\vskip1mm
The next lemma controls $v$ and its derivatives in terms of $\ell$ and the upper barrier $B$. Recall that
$\bar{\rho}_i$ is (up to a factor) the parabolic distance to the singularity at $x^i$.

\begin{Lem}	\label{new:ratios}
There is $c>0$ such that on $\cP_{c\ell}(x,t)$ the flow has $|A|\le1/(c\ell)$. At $(x,t)$, 
\[
|\nabla^Tv|\le C\,\frac B\ell\, ,\qquad |\Hess^T_v|+|\partial_tv|\le C\,\frac B{\ell^2}\, .
\]
Moreover, near each $x^i$,
\[
	\frac Bv\le C\bar\rho_i^{\,-4\theta},\qquad \frac1{\ell v}\le C\bar\rho_i^{\,2\alpha_i-1-2\theta},
\]
and both are bounded elsewhere.
\end{Lem}

\begin{proof}
The set where $\chi=1$ contains   the neighborhood $U'$ of the generic  singularities (Section~\ref{new:sec-upper}); away from these points and the $x^i$'s the curvature is bounded, so $\ell=1$ is the right scale. Where two cases
overlap, the values are comparable. For $c$ small, on $\cP_{c\ell}(x,t)$ we have $|A|\le1/(c\ell)$, and $\ell$, $\bar\rho_i$ and $B$ change by bounded factors, and so does $H$ where $\ell=1/H$, by Lemma~\ref{l:Rss}(i), Lemma~\ref{l:Rss}(iv) on the conical
ends, Lemma~\ref{l:Rfarcyl}(1) on $F^i$, Lemmas~\ref{l:RNSnbhd}(1)--(2) where $\chi=1$, and
bounded geometry elsewhere. 

Let $M$ be the supremum of $v$ over $\cP_{c\ell}(x,t)$. By Proposition~\ref{new:vupper}, and since $B$ changes by bounded factors there and is bounded below on the part below $t_1$, where the flow is smooth, $M\le CB(x,t)$. So Lemma~\ref{new:gradcrit} with $r=c\ell$ gives the bounds on $\nabla^Tv$ and $\Hess^T_v$, and $\partial_tv=\Delta v+|A|^2v$ the bound on $\partial_tv$. The bounds on the ratios follow from Theorem~\ref{new:twosided} and $v\le NB$. The upper and lower bounds of
Theorem~\ref{new:twosided} differ by a factor $C\tau^{-2\theta}$ on the cores and $C\rho_i^{-4\theta}$ on the ends; $v\ge c\tau^{\theta-\alpha_i}$ on the cores, $v\ge c\rho_i^{-2\alpha_i+2\theta}$
on
the conical ends and $v\ge c\rho_i^{1-2\alpha_i+2\theta}H$ on $F^i$; and elsewhere
$B\le C\max\{1,H\}$ and $v\ge\max\{v_0,|H|/(2\kappa)\}$.
\end{proof}

The three parts of the Hessian are
bounded by $(B/v)^k/(\ell v)$ with $k=0,1,2$, as in the Introduction (Theorems~\ref{new:mixed} and~\ref{new:dn}),
which explains the losses $2\theta$, $6\theta$ and $10\theta$ in their decay.

\begin{Thm}	\label{new:mixed}
The mixed part $m$ is bounded on $\Sigma_t$ for $t<T$. Moreover, near each $x^i$,
\[
m\le C\bar\rho_i^{\,2\alpha_i-1-6\theta},\qquad
\frac{|A|}{v}+\frac{|H|}{v}\le C\bar\rho_i^{\,2\alpha_i-1-2\theta}.
\]
\end{Thm}

\begin{proof}
On $[0,t_1]$ the flow is smooth, so let $t\ge t_1$. By Lemma~\ref{new:ratios}, $|\nabla^Tv|\le CB/\ell$ and $|A|\le C/\ell$, so
\[
m=\frac{|\nabla^Tv|}{v^2}\le C\,\frac Bv\cdot\frac1{\ell v},\qquad
\frac{|A|}v+\frac{|H|}v\le\frac C{\ell v},
\]
and the claims follow from the bounds on the ratios in Lemma~\ref{new:ratios}.
\end{proof}

\begin{proof}[Proof of Corollary~\ref{c:Rharnack}]
On $[0,t_1]$ the flow is smooth. For $t\ge t_1$, the first bound is Theorem~\ref{new:mixed}, since $|\nabla^T\log v|=m\,v$. For the second, Lemma~\ref{new:ratios} gives $|\partial_tv|\le CB/\ell^2$, so $|\partial_t\log v|/v^2\le C(B/v)(\ell v)^{-2}$, which is bounded by Lemma~\ref{new:ratios}, since $\theta\le(2\alpha_i-1)/12$.
\end{proof}

\section{The double normal part}	\label{s:Rdn}

The double normal part is $\Hess_u(\nn,\nn)=-v_s/v^3$. We will bound  it by using the maximum principle to bound $v_s$ by a barrier of size $B^2/\ell$.  Unlike $v$, $v_s$ is not a Jacobi field and we will need different barriers.  However, the main part of the equation for $v_s$ is the Jacobi operator and the construction will be similar in spirit.

\subsection{The equation}	\label{ss:Rdneq}

The variation $v_s$ of $v$ across the levels solves an inhomogeneous Jacobi equation with a source.

\begin{Pro}	\label{l:Rslave}\label{p:vs}
Assume $(\star)$. With $\Box=\partial_t-\Delta^T$, for every function $w$ on the levels
\begin{align}	\label{e:Rcomm}
	[\partial_s,\Box]\,w=-2\,A(\nabla^Tv,\nabla^Tw)-2\,v\,\langle A,\Hess^T_w\rangle\, .
\end{align}
Consequently, on $\Sigma_t$ and on every level of $u$ near $\Sigma_t$,
\begin{equation}\label{new:vs-eq}
(\partial_t-\Delta-|A|^2)\,v_s=S\, ,
\end{equation}
where
	$S=2\,v^2\,\Tr_g A^3+4\,v\,\langle A,\Hess^T_v\rangle+2\,A(\nabla^Tv,\nabla^Tv)\, .$
\end{Pro}

\begin{proof}
In normal coordinates at a point, \eqr{e:q1} gives $\partial_sg^{ij}=2vA_{ij}$, and with the Codazzi equation $\partial_s\Gamma^k_{ii}=v\,\nabla_kH-2A(\nabla^Tv)_k-H\,\nabla_kv$. Hence
\begin{align}
	[\partial_s,\Delta^T]\,w=2\,v\,\langle A,\Hess^T_w\rangle+2\,A(\nabla^Tv,\nabla^Tw)+H\,\langle\nabla^Tv,\nabla^Tw\rangle-v\,\langle\nabla^TH,\nabla^Tw\rangle\, .\notag
\end{align}
In Euclidean coordinates $\partial_t=\frac{\partial}{\partial t}-H\nn_i\partial_{x_i}$ along the flow and $\partial_s=v\,\nn_j\partial_{x_j}$. Commuting them, and using $\partial_t\nn=\nabla^TH$, $\partial_s\nn=-\nabla^Tv$ from \eqr{e:q1}, and $\partial_tv+\partial_sH=0$, which is \eqr{e:Jeqn} together with \eqr{e:q1}, gives $[\partial_s,\partial_t]\,w=H\,\langle\nabla^Tv,\nabla^Tw\rangle-v\,\langle\nabla^TH,\nabla^Tw\rangle$. Subtracting gives \eqr{e:Rcomm}. For the second claim, $\Box v_s=\partial_s(\Box v)-[\partial_s,\Box]v$, with $\Box v=|A|^2v$, \eqr{e:q5} and \eqr{e:Rcomm} for $w=v$. This gives $\Box v_s=|A|^2v_s+S$. The proof uses only the Jacobi equation for $v$ and the variational formulas of Section~\ref{ss:Rjacobi}, which hold on every level.
\end{proof}

We call \eqr{new:vs-eq} the Jacobi form of the equation for the double normal
part. The Jacobi operator in it does not involve $v$. Replacing $u$ by $-u$ changes the signs of $\nn$, $A$, $v_s$ and $S$, so \eqr{new:vs-eq} holds for $\pm v_s$ with sources $\pm S$.
Under the hypotheses of Lemma~\ref{new:gradcrit}, moreover,
\begin{align}	\label{new:Sbound}
	|S|\le\frac{CM^2}{r^3}\qquad\text{at }(x,t)\, ,
\end{align}
on $\Sigma_t$ and on every level of $u$ near $\Sigma_t$. Indeed  
$|S|\le2v^2|A|^3+4v|A|\,|\Hess^T_v|+2|A|\,|\nabla^Tv|^2$, and $v\le M$ and $|A|\le1/r$, while Lemma~\ref{new:gradcrit} bounds the derivatives of $v$.

\subsection{The bowl barrier}	\label{ss:Rbarrier}

Near generic singularities,  $H^3$  would be a natural choice for a barrier. This would work on cylindrical regions, but it does not 
quite work near the tip of a bowl.  We deal with this next. 

Given a Jacobi field $w>0$ and $a\in\RR$,  the chain rule gives
\begin{align}	\label{e:Rva}
	(\partial_t-\Delta-|A|^2)\,w^{3+a}=(a+2)\Big(|A|^2-(a+3)\,|\nabla\log w|^2\Big)\,w^{3+a}\, .
\end{align}
With $w=H$ and $a=0$, $H^3$ gains $2(|A|^2-3|\nabla\log H|^2)$ times itself, which is negative on a shoulder around the tip of a bowl.

On a translating soliton with velocity $Ve$, $|e|=1$, oriented so that $H>0$, we have $H=-V\langle\nn,e\rangle$, so $\nabla H=VA(e^T)$, and $|e^T|^2+\langle\nn,e\rangle^2=1$ gives $V^2=H^2+|A^{-1}\nabla H|^2$ wherever $A$ is invertible. On a surface with $H>0$ and Gauss curvature $K>0$, define $\mathcal V\geq H$ by $\mathcal V^2=H^2+|A^{-1}\nabla H|^2$. Then $\mathcal V=V$ on a translating bowl, and $\mathcal V=H$ on a round sphere. Near a cylinder $A$ degenerates, so we will cut off on the cylindrical region. 

{\it The bowl.} 
Let $R_{\max}< 1$ be given by Lemma~\ref{l:Rbowlineq}.
Fix $a\in(-2,0)$ with $(a+3)R_{\max}<1$, and set $c_a=\frac12(a+2)\big(1-(a+3)R_{\max}\big)>0$. On the bowl, $V/H=(1+(rp)^2)^{1/2}$ and $K/H^2=p(1-p)$, with $p$ from  Lemma~\ref{l:Rbowlineq}.  That lemma gives that the first increases 
(since $p \in (1/2,1)$ is increasing) 
and the second decreases with $r$. By \eqr{e:RRbowl} and Lemma~\ref{l:Rbowlineq}(2), $R\leq D^2\to0$ as $r\to\infty$, so we may fix $r_M$ with $R\leq\frac16$ for $r\geq r_M-\frac12$. Let $M_B$ be the value of $(V/H)^{-a}$ at $r=r_M$, and $k_B$ the value of $K/H^2$ at $r=r_M+\frac12$, on the bowl. 

{\it The barrier.}
We will next define a function $F$ that is constant away from the tip region and behaves like $H^a$ near the tip.
Write $q=K/H^2$ and $z=(\mathcal V/H)^{-a}$, and fix a smooth $\beta:\RR\to[0,1]$ with $\beta=0$ on $(-\infty,k_B/2]$ and $\beta=1$ on $[k_B,\infty)$. Define $F$ and $h$ by 
\begin{align}	\label{e:Rh}
	 F=M_B+\beta(q)\,(z-M_B) {\text{ and }} h=\min\{M_B,F\}.
\end{align}
Then $F$ is smooth where $H>0$, since $z$ enters only where $K\ge\frac{k_B}2H^2>0$, so that $A$ is invertible; under (MC), moreover, $|A^{-1}|\le C/H$ there. Since $z\ge1$, $1\le h\le M_B$. On the bowl, $q$ decreases and $z$ increases with $r$, $z=M_B$ at $r=r_M$, and $q=k_B$ at $r=r_M+\frac12$. So $z>M_B$, where $k_B/2\le q\le k_B$, and $h=M_B$ there. Hence on the bowl $H^3h=H^{3+a}V^{-a}$ for $r\le r_M$ and $H^3h=M_BH^3$ for $r\ge r_M$.

\begin{Lem}	\label{l:Rsuper}
There is $\varepsilon'>0$ so that:  At every point at which the flow, rescaled by $H$, is $\varepsilon'$-close in $C^5$, on the backward parabolic ball of radius $1/\varepsilon'$, to a round shrinking sphere or cylinder or a translating bowl, $h$ of \eqr{e:Rh} is defined and
\[
(\partial_t-\Delta-|A|^2)(H^3h)\ge\frac{c_a}{2}\,H^2\cdot H^3h ,
\]
in the viscosity sense.
\end{Lem}

\begin{proof}
Near the models, $h=\min\{M_B,(\mathcal V/H)^{-a}\}$ where $q\ge k_B$, and $h=M_B$ where $q<k_B$, since $z>M_B$ near the part of the bowl where $k_B/2\le q\le k_B$. So $H^3h$ is the minimum of $M_BH^3$ and $H^{3+a}\mathcal V^{-a}$ where $q\ge k_B$, and $M_BH^3$ elsewhere. The minimum of two supersolutions is a supersolution, so it suffices to check the inequality for each case. Both sides of the inequality, divided by $H^5h$, are continuous functions of the rescaled $5$-jet of the flow, so it suffices to check them on the models, with room.
\begin{itemize}
\item $M_BH^3$ is active on cylinders, and on the bowl for $r\geq r_M$. By \eqr{e:Rva} with $a=0$ it gains $2(|A|^2-3|\nabla\log H|^2)\geq|A|^2$ times itself, since $\nabla H=0$ on cylinders and $R\leq\frac16$ on the bowl for $r\geq r_M-\frac12$.
\item $H^{3+a}\mathcal V^{-a}$ is active on spheres, where $\mathcal V=H$, so that it is $H^3$ and gains $2|A|^2$ times itself; and on the bowl for $r\leq r_M$, where $\mathcal V=V$ is constant, so that by \eqr{e:Rva} and Lemma~\ref{l:Rbowlineq}(2) it gains $(a+2)\big(|A|^2-(a+3)|\nabla\log H|^2\big)\geq2c_a|A|^2$ times itself.
\end{itemize}
Since $|A|^2\geq\frac12H^2$ and $c_a\leq1$, this gives the inequality for $\varepsilon'$ small.
\end{proof}

\subsection{The barrier for $v_s$}\label{new:sec-dn}
Each piece of our barrier will
correspond to a piece of the barrier $B$ for $v$. Let $\chi$, $\hat P$ and $\hat\Psi$ be as in
Section~\ref{new:sec-upper}. Let $\rho_0$ be given by Lemma~\ref{l:RNSnbhd}(2) with $\epsilon_0=\varepsilon'$, and let $h$ be as in \eqr{e:Rh}; it is defined wherever $H>0$, and $1\le h\le M_B$. Away from the cores let
\[
B_{V,e}=e^{\Lambda(t-t_1)}\,\hat\Psi^2\times
\begin{cases}
\hat P\,\psi_1(\rho_i)/(1+\eta) & \text{on the conical ends of } x^i,\\
H^3h & \text{on } F^i\cap\{\rho_i\le3r_1\},\\
\chi H^3h+1-\chi & \text{elsewhere.}
\end{cases}
\]
These agree where the cases meet, since $\chi=1$ on
$F^i\cap\{2r_1\le\rho_i\le3r_1\}$, $\chi=0$ on the conical ends, and $\hat P=1+\eta$ and
$\hat\Psi=\psi_1=1$ for $\rho_i\ge2r_1$. On the cores $\{|y_i|\le2R^i_c\}$, and on the
compact components $\mathcal C_t$, let
$B_{V,c}=K^i_2\tau^{-\alpha_i-\frac12-2\theta}\Phi_i$, where $K^i_2$ is fixed below by Lemma~\ref{new:glue} when $S^i$ is noncompact, and $K^i_2=1$ on $\mathcal C_t$.

Since $\psi_a\psi_b=\psi_{a+b}$, comparing with the pieces of $B$ in
Section~\ref{new:sec-upper} gives
\begin{align}	\label{new:square}
	c\,B^2/\ell\le B_{V,e}\le C\,B^2/\ell\quad\text{and}\quad c\,B^2/\ell\le B_{V,c}\le C\,B^2/\ell
	\quad\text{where they are defined}.
\end{align}
On the cores both sides are $\tau^{-2\alpha_i-\frac12-2\theta}$ times functions of $y_i$
that are bounded above and below. On the conical ends with $\rho_i\le r_1$ both are
comparable to $\rho_i^{-4\alpha_i-1-4\theta}$, on $F^i\cap\{\rho_i\le r_1\}$ to
$H^3\rho_i^{2-4\alpha_i-4\theta}$, where $\chi=1$ to $H^3$, and elsewhere to $1$.

\begin{Lem}	\label{new:vslemma}
If $\Lambda$ is large, depending on $r_1$ and on the flow away from the $x^i$, then
$B_{V,e}$ and $B_{V,c}$ are positive strict supersolutions of the Jacobi equation where
they are defined, for $t_1\le t<T$: for a constant $c>0$ depending on $r_1$,
\[
(\partial_t-\Delta-|A|^2)B_{V,e}\ge c\,\ell^{-2}B_{V,e},\qquad
(\partial_t-\Delta-|A|^2)B_{V,c}\ge c\,\ell^{-2}B_{V,c}.
\]
\end{Lem}

\begin{proof}
{\it Cores and compact components.} This is Lemma~\ref{new:core} with $\gamma=\alpha_i+\frac12+2\theta$, which applies on the cores and on $\mathcal C_t$ by Section~\ref{app:constants}(a). It gives the gain $\frac\gamma{2\tau}B_{V,c}$; here $\ell=\sqrt\tau$.

{\it Conical ends, $\rho_i\le r_1$.} Here $\hat P=P_+$ and $\hat\Psi^2\psi_1=a\rho_i^{-p}$,
with $p=4\alpha_i+1+4\theta$ and a constant $a>0$. So this is Lemma~\ref{new:radial}(2)
with this $p$, whose hypotheses hold by Lemma~\ref{l:Rss}(iv) and since $K_PR_0\ge4+6(p^2+C_g^2)$, as $R_0$ is large. Its gain is
$\frac14\eta R_0|y_i|\rho_i^{-p-2}\ge\frac14K_PR_0\,\rho_i^{-2}\rho_i^{-p}$, up to the
factor $a/(1+\eta)$; here $\ell=\rho_i$.

{\it The far parts, $\rho_i\le3r_1$.}$^\dagger$ On $F^i\cap\{\rho_i\le3r_1\}$ the hypothesis of
Lemma~\ref{l:Rsuper} holds by Lemma~\ref{l:Rfarcyl}(1). So
$(\partial_t-\Delta-|A|^2)(H^3h)\ge\frac{c_a}{2}H^2\cdot H^3h$. Moreover
$\hat\Psi^2=\psi_{q_f}$, with $q_f=4\alpha_i-2+4\theta$, so
$\rho_i^2|(\partial_t-\Delta)\hat\Psi^2|\le C(q_f)\hat\Psi^2$ and
$|\nabla\hat\Psi^2|\le q_f\hat\Psi^2/\rho_i$; and
$|\nabla(H^3h)|\le(3C_A+C_h)H\cdot H^3h$, where $C_h$, which depends only on the models, bounds $|\nabla\log h|/H$ near them. Since $H\rho_i\ge cR_0$, and $(cR_0)^2c_a\ge32\,C(q_f)$ and $cR_0c_a\ge32(3C_A+C_h)q_f$ as $R_0$ is large, the product rule
gives
$(\partial_t-\Delta-|A|^2)(H^3h\,\hat\Psi^2)\ge\frac{c_a}{4}H^2\cdot H^3h\,\hat\Psi^2$;
here $\ell=1/H$.

{\it Elsewhere.} Where $\chi=1$ we are in the mean convex neighborhoods $\mathcal U$ of Lemma~\ref{l:RNSnbhd}, where $B_{V,e}=e^{\Lambda(t-t_1)}H^3h$ and $\ell=1/H$. There are three cases.
\begin{itemize}
\item {\it $\chi=1$ and $H\ge1/\rho_0$.} The hypothesis of Lemma~\ref{l:Rsuper} holds by Lemma~\ref{l:RNSnbhd}(2); this gives the gain $\frac{c_a}{2}\ell^{-2}B_{V,e}$.
\item {\it $\chi=1$ and $H\le1/\rho_0$.} The geometry is bounded and $H\ge h_0$ by (MC). The function $h$ is smooth where $H>0$, apart from the minimum, and $(\partial_t-\Delta)\mathcal V$ involves only $A$, $\nabla A$ and $\nabla^2A$: the third derivatives of $A$ in $\partial_t\mathcal V$ and $\Delta\mathcal V$ cancel, since $[\partial_t-\Delta,\nabla]$ is of lower order. So (MC) gives $(\partial_t-\Delta-|A|^2)(H^3h)\ge-C\,H^3h$ for each branch, and the factor $e^{\Lambda(t-t_1)}$ absorbs this, together with the gain $\ell^{-2}B_{V,e}\le\rho_0^{-2}B_{V,e}$.
\item {\it $\chi<1$, and where $\hat P$, $\hat\Psi$ and $\psi_1$ vary on the conical ends.} Where $\chi=0$, $(\partial_t-\Delta-|A|^2)e^{\Lambda(t-t_1)}=(\Lambda-|A|^2)e^{\Lambda(t-t_1)}$. Where $\chi$ varies, and where $\hat P$, $\hat\Psi$ and $\psi_1$ vary, all other terms are bounded below by $-C\,B_{V,e}$, and $\Lambda$ absorbs them (on the conical ends the term $\partial_t\hat P\ge0$ is unbounded as $t\to T$, but has the favorable sign); there $\ell$ is bounded below by a constant depending on $r_1$.
\end{itemize}
\end{proof}

At fixed $|y_i|$ all pieces carry the power $\tau^{-2\alpha_i-\frac12-2\theta}$, as the
comparison after \eqr{new:square} shows. So Lemma~\ref{new:glue} applies, with
$B_c=\tau^{-\alpha_i-\frac12-2\theta}\Phi_i$, $B_e=B_{V,e}$, and its constant $\kappa$ as $K^i_2$; here
$a-\gamma=-(2\alpha_i+1+4\theta)$ on the conical ends and
$a-\gamma=-(2\alpha_i-1+4\theta)$ on the cylindrical ends, where $h=M_B$ on the
graphical part. Let $B_V$ be $K$ times the glued
barrier: $B_{V,c}$ on $\{|y_i|\le R_0\}$ and on $\mathcal C_t$,
$\min\{B_{V,c},B_{V,e}\}$ on $\{R_0\le|y_i|\le2R^i_c\}$, and $B_{V,e}$ elsewhere.  

We are now ready to prove the double normal bound.  It is useful for the sequel to separate out the maximum principle step:

\begin{Pro}	\label{new:vs}
If $K$ is large, depending on the flow at $t_1$, then
$(\partial_t-\Delta-|A|^2)B_V>|S|$ on $\Sigma_t$ for $t_1\le t<T$, and $|v_s|\le B_V$ for $t_1\le t<T$.
\end{Pro}

\begin{proof}
Since $|S|\le2v^2|A|^3+4v|A|\,|\Hess^T_v|+2|A|\,|\nabla^Tv|^2$, Lemma~\ref{new:ratios} and $v\le NB$ give $|S|\le CB^2/\ell^3$. By \eqr{new:square} this is at most a fixed
multiple of $\ell^{-2}B_{V,e}$, respectively $\ell^{-2}B_{V,c}$, and so of the gain in
Lemma~\ref{new:vslemma}. So $(\partial_t-\Delta-|A|^2)B_V>|S|$ for $K$ large. Take $K$
also so large that $|v_s|<B_V$ at $t=t_1$. Since $(\partial_t-\Delta-|A|^2)(\pm v_s)=\pm S$, the first-touching argument of Section~\ref{ss:Rtools}, applied to $\pm v_s$ and $B_V$, gives $|v_s|<B_V$.
\end{proof}

\begin{Thm}	\label{new:dn}
{\it The double normal part.}
$\Hess_u(\nn,\nn)$ is bounded on $\Sigma_t$ for $t<T$, and
$|\Hess_u(\nn,\nn)|\le C\bar\rho_i^{\,2\alpha_i-1-10\theta}$ near each $x^i$.
\end{Thm}

\begin{proof}
By Proposition~\ref{new:vs} and \eqr{new:square},
\[
|\Hess_u(\nn,\nn)|=\frac{|v_s|}{v^3}\le\frac{B_V}{v^3}\le C\Big(\frac Bv\Big)^2\frac1{\ell v}.
\]
The claims follow from Lemma~\ref{new:ratios}, since $\theta\le(2\alpha_i-1)/12$.
\end{proof}

\subsection{The starting time}	\label{app:constants}

All constants of Sections~\ref{s:Rmixed} and~\ref{s:Rdn} except the starting time $t_1$ and the overall factors of the barriers are fixed above, each depending only on the $S^i$, $\theta$ and the constants fixed before it. None of them uses the flow at $t_1$, and no barrier depends on $v$. The starting time is fixed last, by smallness conditions collected here.

\smallskip
\noindent(a) Take $s_1$ so large that, with $\tau_1=e^{-s_1}$:
\begin{itemize}
\item $\varepsilon\le\varepsilon_C(2R^i_c,\theta)$ on $[s_1-1,\infty)$
(Lemma~\ref{new:core}) at each noncompact $S^i$, and
$\varepsilon\le\varepsilon_C(R_{S^i},\theta)$ on $[s_1-1,\infty)$ at each compact $S^i$, where $R_{S^i}\ge1$ and $S^i\subset B_{R_{S^i}}$, as in the proof of Lemma~\ref{new:core}; since $\varepsilon_C(R,\gamma)$ is increasing in $\gamma$ (the remark after Lemma~\ref{new:core}), these conditions also hold with $\alpha_i+\frac12+2\theta$ in place of $\theta$, which Lemma~\ref{new:vslemma} uses;
\item $\hat R(s_1)\ge4R^i_c$ and $2R^i_c\sqrt{\tau_1}\le r_1$;
\item $\tau_1\le\min\{\delta_\Upsilon,1/C_4,1/\Lambda\}$
and $6K_P\sqrt{\tau_1}\le r_1$;
\item $\tau_1$ is below the $\delta$ of Lemma~\ref{l:Rcompact} at each compact $S^i$.
\end{itemize}
Set $t_1=T-\tau_1$.

\smallskip
\noindent(b) The overall factors of Propositions~\ref{new:vupper}, \ref{new:vlower}
and~\ref{new:vs}, from the flow at $t_1$.

\begin{proof}[Proof of Theorem~\ref{t:Rintro}]	 
 On $[0,t_1]$ the flow is smooth,
and $u$ is smooth near $\Sigma_t$ by Lemma~\ref{l:Rsetup}, so $\Hess_u$ is
bounded there. For $t_1\le t<T$, the tangential part is bounded by
Proposition~\ref{p:Rtang}, the mixed part by Theorem~\ref{new:mixed}, and the double
normal part by Theorem~\ref{new:dn}. 
\end{proof}

\appendix

\section{Polynomial bounds and separation of the ends}	\label{new:appA}\label{s:Rapp}

Lemma~\ref{new:polyb} says that on each end the first eigenfunction $\phi_1$ grows like the power $|y|^\gamma$ of the end, up to factors $|y|^{\pm\delta}$. We prove it by comparing $\phi_1$ with $|y|^{\gamma\pm\delta}$, which are sub- and supersolutions far out on the end.

\begin{proof}[Proof of Lemma~\ref{new:polyb}]
Let $E$ be an end of $S$. On a shrinker $\langle y,\nn\rangle=2H$ and $\Delta y=-H\nn$,
so $\Delta|y|^2=4-4H^2$ and $|y^T|^2=|y|^2-4H^2$. For $p\in\RR$ this gives
\[
L|y|^p=\Bigl(\frac12-\frac p2+|A_S|^2\Bigr)|y|^p+O(|y|^{p-2}),
\]
with $O$ depending on $p$. On a conical end, $|A_S|^2\to0$, and on a
cylindrical end $|A_S|^2\to\frac12$. Since
$\alpha=-\mu_1-\frac12$, this gives
\[
(L+\mu_1)|y|^{\gamma\pm\delta}=\Bigl(\mp\frac\delta2+o(1)\Bigr)|y|^{\gamma\pm\delta}
\qquad\text{as }|y|\to\infty\text{ on }E .
\]
Take $R_\delta$ so large that on $E_\delta=E\cap\{|y|\ge R_\delta\}$ the factor
$\mp\frac\delta2+o(1)$ has the sign of $\mp1$, that
$|A_S|^2+\frac12+\mu_1<0$ (recall $\mu_1<-1$ by \cite{CM7}), and that 
$\partial E_\delta$ is $E\cap\{|y|=R_\delta\}$. Let $f=\phi_1-C_\delta|y|^{\gamma+\delta}$, or
$f=c_\delta|y|^{\gamma-\delta}-\phi_1$, with the constant chosen so that $f\le0$ on
$\partial E_\delta$. Then $(L+\mu_1)f\ge0$ on $E_\delta$, since
$(L+\mu_1)\phi_1=0$. Both $\phi_1$ and $|y|^{\gamma\pm\delta}$ lie in $W^{1,2}(E_\delta,e^{-|y|^2/4})$. Multiply by $f^+e^{-|y|^2/4}$ and integrate by
parts over $E_\delta$; this is allowed since $f^+=0$ on $\partial E_\delta$ and the
Gaussian weight dominates the polynomial growth. This gives
\[
-\int_{E_\delta}|\nabla f^+|^2e^{-|y|^2/4}
+\int_{E_\delta}\Bigl(|A_S|^2+\frac12+\mu_1\Bigr)(f^+)^2e^{-|y|^2/4}\ge0 .
\]
Both terms are nonpositive, and the second is negative unless $f^+\equiv0$. So
$f\le0$ on $E_\delta$ in either case, which gives the claim.
\end{proof}

\vskip1mm
We will also need the following consequence of \cite{BK}, cf. \cite{Wa}.

\begin{Cor}	\label{c:BK1}
{\it Separation of the ends.} There is $c_S>0$, depending only on $S$, such that for $|y|$
large the ball of radius $c_S|y|$ about a point $y$ of an end of $S$ meets $S$ only in
that end.
\end{Cor}

\begin{proof}
By \cite[Theorem~1.10]{BK}, $S$ has finitely many ends, each within $o(|y|)$ of the cone over a set $\Gamma_E\subset S^2$ (for a cylindrical end, $\Gamma_E=\{v_E\}$ and the cone is the ray $\RR_+v_E$), and the $\Gamma_E$ are pairwise disjoint compact subsets of $S^2$. So they are at positive distance from each other, and this gives $c_S$.
\end{proof}

\section{The local geometry}	\label{app:geom}

We prove here the lemmas of Sections~\ref{s:RNS}--\ref{new:sec-radial} that describe the flow near the singular points, and the finiteness of the non-generic singular points (Lemma~\ref{new:finite}). 

\begin{Lem}	\label{new:finite}
 There are finitely many non-generic singular points at time $T$.
\end{Lem}

\begin{proof}
The flow on $[0,T]$ is almost regular \cite[Theorem~1.1]{BK}. By \cite[Theorem~1.9(b),(c)]{BK},
its non-generic singular points form a closed set, and each is backward isolated, so in
particular isolated in the slice $\{t=T\}$. Since the singular set at time $T$ is compact, there
are finitely many.
\end{proof}

The next lemma is a quantitative form of the mean convex neighborhood theorem that is  implicit in \cite{ChHH}.

\begin{Lem}	\label{l:RCHHq}
For every $\epsilon>0$, $L_Q\geq2\sqrt2$ and $\Lambda_0<\infty$ there are $R_Q$ and $\delta>0$ with the following property. Let $M_t$ be a mean curvature flow of closed embedded surfaces with entropy at most $\Lambda_0$, and suppose that on $B_{R_Q\sqrt{2\tau}}(y_0)\times[t_*-R_Q^2\tau,t_*]$ the flow is, after rescaling by $\sqrt{\tau}$ about $(y_0,t_*)$, $\delta$-close in $C^2$, as a normal graph of multiplicity one, to the round cylinder shrinking to its axis at time $t_*+\tau$. Then at every point of the flow in $B_{L_Q\sqrt{\tau}}(y_0)\times[t_*,t_*+\tau)$ the mean curvature is positive, and the flow rescaled by $|A|$ about the point is $\epsilon$-close in $C^{[1/\epsilon]}$, on the backward parabolic ball of radius $1/\epsilon$, to a round shrinking sphere or cylinder, a translating bowl, or an ancient oval.
\end{Lem}

\begin{proof}[Sketch of proof]
This follows from the argument in \cite[Section~8]{ChHH}, applied to a sequence of flows. If the lemma fails, there are flows $M^i_t$ as in the statement with $\delta_i\to0$ and $R_{Q,i}\to\infty$, and points $(p_i,s_i)$ where the conclusion fails.
Rescale so that the model is the cylinder of radius $\sqrt2$ that pinches at time $1$. By Huisken's monotonicity and the entropy bound $\Lambda_0$, every limit of the $M^i_t$, and of their rescalings about $(p_i,s_i)$ by the regularity scale $r_i$, is an ancient, unit-regular, cyclic integral Brakke flow with entropy at most that of the cylinder.
The part of the flow outside the region of closeness contributes at most $C\Lambda_0\e^{-R_{Q,i}^2/C}$. By \cite[Theorem~1.2]{ChHH}, the first limit, which is the model for $t\leq0$, is the model up to time $1$, so by \cite{W1} the flows converge smoothly on compact subsets of $\RR^3\times(-\infty,1)$. Hence $s_i\to1$ and $r_i\to0$. The second limit is then either a plane, which contradicts the definition of $r_i$, or a round sphere or cylinder, a bowl or an ancient oval, which gives the conclusion at $(p_i,s_i)$ for $i$ large.
\end{proof}

\begin{proof}[Proof of Lemma~\ref{l:RNSnbhd}]
Each generic singularity has a mean convex neighborhood \cite[Theorem~1.7]{ChHH}, so sets $U$ as in the lemma exist. Since $T$ is the first singular time, \cite[Theorem~1.6]{ChHH} also shows that every limit flow at a point of $\mathcal S_T\cap\overline U$ is a round shrinking sphere or cylinder or a translating bowl.

A sphere occurs only at a sphere point. Indeed, a sphere limit along points $(p_i,t_i)\to(z,T)$ makes the component of $\Sigma_{t_i}$ through $p_i$ a closed convex surface of diameter at most $C/H(p_i,t_i)$. Since $\Sigma_t$ has finitely many components for $t<T$, one of them then shrinks to $z$ at time $T$, to a round point by \cite{H1}. So the tangent flow at $z$ contains a round sphere, and $z$ is a sphere point.

The regularity scale tends to $0$ at points approaching $(\mathcal S_T\cap\overline U)\times\{T\}$, and $\mathcal S_T\cap\overline U$ is compact. So the usual argument by contradiction, rescaling by the regularity scale and using \cite{W1}, gives for each $\epsilon_0>0$ the closeness of (2), even with $|A|$ in place of $H$, at all points of $\mathcal U$ in a space-time neighborhood $\mathcal N=\mathcal N(\epsilon_0)$ of $(\mathcal S_T\cap\overline U)\times\{T\}$. Since the curvature is bounded on $\mathcal U$ outside $\mathcal N$, every point of $\mathcal U$ with $H\geq1/\rho_0$ lies in $\mathcal N$ for $\rho_0$ small, and (2) follows.

Fix $\mathcal N$ for one small $\epsilon_0$. On $\mathcal N$, as on the models, $H>0$ and $H$ is large, $|A|\leq CH$, and $|\nabla^kA|\leq CH^{k+1}$ for $k\leq2$. Outside $\mathcal N$ the flow is smooth up to $T$, with bounded geometry. There $H>0$, also at time $T$, by strict nesting in the mean convex neighborhoods and the strong maximum principle. So $H\geq h_0>0$ there by compactness, and (1) follows.
\end{proof}

\begin{proof}[Proof of Lemma~\ref{l:Rss}]
{\it (i).} Since $S$ is unique, the rescaled flow converges to $S$, smoothly on compact sets since $S$ is smooth and has multiplicity one \cite{W1}. So, by a diagonal argument, the rescaled flow is $C^2$-close to $S$ with an error $\varepsilon_0(s)\to0$ on balls of radius $R(s)\to\infty$ (cf.\ \cite{CM1} for how such graphical regions expand), and (i) holds with $\varepsilon=\max\{\varepsilon_0(s),1/\log R(s)\}$.

{\it (iii).} We may replace $\varepsilon$ by a larger function. We do this to make $\varepsilon^{-1}$ $\frac14$-Lipschitz.
Replace $1/\varepsilon(s)$ by $\inf_{s'\leq s}\big(1/\varepsilon(s')+(s-s')/4\big)$; this makes $\varepsilon$ larger and still tends to $0$. Then $2^{-1/2}\min\{\e^{1/\varepsilon},\e^{s/4}\}\leq\hat R\leq\min\{\e^{1/\varepsilon},\e^{s/4}\}$. Since $1/\varepsilon$ is nondecreasing and $\frac14$-Lipschitz, so is $\log\hat R$, as a smooth minimum of $1/\varepsilon$ and $s/4$.
So the static radius $\hat R(s)\e^{-s/2}$ is continuous and strictly decreasing, and its inverse, the entry time $s_\rho$, is continuous and decreasing in $\rho$. Finally the bound for $\varepsilon$ in (iii) holds, since for $s-2\log\hat R(s)\leq s'\leq s$ we have $1/\varepsilon(s)-1/\varepsilon(s')\leq(s-s')/4\leq\frac12\log\hat R(s)\leq1/(2\varepsilon(s))$.

{\it (ii).} Let $\mathcal C$ be the asymptotic cone of a conical end of $S$. Since the end is smoothly asymptotic to $\mathcal C$ \cite[Theorem~1.10]{BK}, there is a nonincreasing $\eta_S(\Xi)\to0$ such that, for $\Xi$ large, the part of the end in $\{\Xi/2\leq|y|\leq2\Xi\}$ is the normal graph over $\mathcal C$ of a function $w_S$ with $|w_S|+\Xi|\nabla w_S|+\Xi^2|\nabla^2w_S|\leq\eta_S(\Xi)\,\Xi$. We say that a normal graph $w$ over $x_0+\mathcal C$ near a point at distance $\rho$ from $x_0$ has relative norm at most $\delta$ if $|w|+\rho|\nabla w|+\rho^2|\nabla^2w|\leq\delta\rho$. Fix $x\in x_0+\mathcal C$ with $|x-x_0|=\rho$, let $\hat t_\rho=T-\varepsilon(s_\rho)\rho^2$, and fix constants $K$ large and then $c>0$ small, depending only on $S$ and on the constant $\epsilon_W$ of White's local regularity theorem \cite{W1}. Recall that $s_\rho$ is as large as we like after decreasing $r_0$.

We use a pseudolocality argument, cf.\ \cite{CIM} for shrinkers and \cite{CM1} for the rescaled flow, in the form of \cite[Lemma~9.1]{ChS}; see also \cite[Corollary~1.2]{ChS} and \cite[Section~4]{LZ}.    The times $t\leq\hat t_\rho$ with $T-t\leq4c^2\rho^2$ precede $t_\rho$, and at them (i), (iii), the asymptotics of the end and Corollary~\ref{c:BK1} make the flow on $B_{2Kc\rho}(x)$ a normal graph over $x_0+\mathcal C$ of relative norm at most $C\big(\varepsilon(s_\rho)^{1/2}+\eta_S(1/(2c))\big)$; at time $\hat t_\rho$ it is at most $\delta_1=C\big(\varepsilon(s_\rho)^{1/2}+\eta_S(\varepsilon(s_\rho)^{-1/2})\big)$. By Huisken's monotonicity, Gaussian density ratios at scales at most $c\rho$ about points of $B_{c\rho}(x)\times[\hat t_\rho-c^2\rho^2,T)$ are bounded by ratios over these time slices, where the flow is close to a plane at the scale $\rho$. So, for $K$ large, then $c$ small and $s_\rho$ large, White's theorem and Lemma~\ref{l:EH2}, rescaled, give
\[
(\ast)\qquad |\nabla^kA|\le C\rho^{-1-k}\ \ (k\le2)\quad\text{on } B_{3c\rho/8}(x)\ \text{ for } t\ge\hat t_\rho-c^2\rho^2/8,
\]
with a single sheet there. Since Huisken's monotonicity holds for the whole flow, these bounds on Gaussian density ratios hold at every point of the flow in $B_{c\rho}(x)\times[\hat t_\rho-c^2\rho^2,T)$, however close to $T$. A second sheet within distance $d\leq c\rho$ at such a time would give a Gaussian density ratio at least $1+c_W$ at the scale $d$, with $c_W>0$ depending only on the curvature bound of \cite{W1}; so, with the bounds on the density ratios below $1+c_W$, there is a single sheet. The interval $[\hat t_\rho,T)$ has length $\varepsilon(s_\rho)\rho^2$, so the flow changes the graph by at most $C\varepsilon(s_\rho)$ in relative norm over it. So for $t\geq t_\rho$ the flow on $B_{c\rho/4}(x)$ is a normal graph over $x_0+\mathcal C$ of relative norm at most $\delta(\rho)=\delta_1+C\varepsilon(s_\rho)$. As $s_\rho$ decreases in $\rho$ and $\varepsilon$ and $\eta_S$ are nonincreasing, $\delta(\rho)$ is nondecreasing, and it tends to $0$ as $\rho\to0$. Covering the part of $x_0+\mathcal C$ at distance $\rho$ from $x_0$ by such balls gives (ii).

{\it (iv).} Let $(x,t)$ be as in (iv), let $x'$ be the point of $x_0+\mathcal C$ nearest to $x$, where $\mathcal C$ is the asymptotic cone of the end, and let $\rho'=|x'-x_0|$. By (ii) on the static part, and by (i) and the asymptotics of the end 
 on the graphical part, $|x-x'|\le c\rho/8$ once $r_0$ is small and $R_G$ is large, so $\rho\le2\rho'$. Take $R_G\ge8/c$ and $c_g\le c/16$. Then $B_{c_g\rho}(x)\subset B_{3c\rho'/8}(x')$, and every $t'\in(t-c_g^2\rho^2,t]$ has $T-t'\le\rho^2(R_G^{-2}+c_g^2)\le c^2\rho'^2/8$, so $t'\ge\hat t_{\rho'}-c^2\rho'^2/8$. So $(\ast)$ at $x'$ gives (iv).

{\it (v).} Since $S$ is connected, for $s$ large the graph in (i) is connected, so it lies on one component of $\Sigma_t$; since its rescalings contain the graph in (i), this component has Gaussian density at least that of $S$ at $(x_0,T)$. The Gaussian density of the flow at $(x_0,T)$ equals that of $S$, since the tangent flow is $S$ with multiplicity one, and it is additive over components. The components of $\Sigma_t$, $t<T$, correspond along the flow and are finitely many. If (v) failed, there would be $t_i\to T$ and points $x_i\to x_0$ of other components of $\Sigma_{t_i}$, which we may take on one component $\mathcal C'$. Then $(x_0,T)$ lies in the closure of the track of $\mathcal C'$, however $T-t_i$ compares with $|x_i-x_0|^2$, so by upper semicontinuity $\mathcal C'$ has Gaussian density at least $1$ at $(x_0,T)$. This is a contradiction.
\end{proof}

\begin{proof}[Proof of Lemma~\ref{l:Rsheets}]
Let $\hat R=\hat R(s_{2\rho})$, so that $\sqrt{T-t_{2\rho}}=2\rho/\hat R$.

By Lemma~\ref{l:Rss}(i), $\Sigma_{t_{2\rho}}\cap B_{2\rho}(x_0)$ is the normal graph over $x_0+\sqrt{T-t_{2\rho}}\,S$ of a function of rescaled $C^2$ norm at most $\varepsilon(s_{2\rho})$. On $\{\rho/4\leq|x-x_0|\leq2\rho\}$, $\sqrt{T-t_{2\rho}}\,S$ lies within $\eta_S(\hat R/8)\rho$ of the cones over its conical ends and within $C\rho/\hat R$ of the axes of its cylindrical ends. 
By \cite[Theorem~1.10]{BK}, far out a cylindrical end of $S$ is a normal graph over the cylinder of radius $\sqrt2$ about its axis, so it lies within $2$ of the axis. So let $x$ with $\rho/2\leq|x-x_0|\leq\rho$ be at distance at least $C\rho/\hat R$ from the axes and at least $c\rho/10$ from the cones, and let $d$ be its distance from their union. Then $d\geq C\rho/\hat R$ and $B_{d/2}(x)\cap\Sigma_{t_{2\rho}}=\emptyset$, for $\hat R$ large. Here $c$ is the constant of the proof of Lemma~\ref{l:Rss}(ii).

By avoidance with the sphere of radius $d/2$ about $x$, 
  $B_{d/4}(x)$ does not intersect $\Sigma_t$ for $t_{2\rho}\leq t<T$. This uses $T-t_{2\rho}=4\rho^2/\hat R^2\leq d^2/32$, which holds for $d\geq C\rho/\hat R$.

The points within $c\rho/10$ of $x_0+\mathcal C$ lie in the balls of the proof of Lemma~\ref{l:Rss}(ii). There $\Sigma_t$ has no other sheets and is the normal graph of Lemma~\ref{l:Rss}(ii).

\end{proof}

\begin{proof}[Proof of Lemma~\ref{l:Rcompact}]
By Remark~\ref{p:Rstruct}(1), $S$ has multiplicity one, so by \cite{Sch} it is the unique tangent flow at $(x_0,T)$ and Lemma~\ref{l:Rss} applies. Since $S$ is compact and $\hat R(s)\to\infty$, for $s$ large the graph in Lemma~\ref{l:Rss}(i) is a normal graph over all of $S$. It is a closed connected embedded surface, so it is a component $\tilde{\mathcal C}_s$ of $\tilde\Sigma_s$. 

 It remains to show that there is some $\delta > 0$ so that the other components avoid $B_\delta(x_0)$ for $T-\delta\leq t<T$.  
 This follows from the avoidance principle using the following fact:
 \begin{itemize}
 \item Since $S$ is compact and $\hat{R}(s) \to \infty$, given any large $N$, we can choose a $\tau_N > 0$ so that
 $B_{N\, \sqrt{\tau_N}}(x_0) \cap \Sigma_{T-\tau_N}$ consists of just one component, this component is a graph over (a scaling of) $S$, and this component is inside of $B_{\frac{N}{4}\, \sqrt{\tau_N}}(x_0)$.
 \end{itemize}
 Since $\partial B_{\frac{N}{2}\, \sqrt{\tau_N}}(x_0) $ avoids $\Sigma_{T-\tau_N}$, they stay disjoint when we flow both simultaneously. For $N$ sufficiently large, the sphere still has positive radius at time $T$. This gives the claim.
\end{proof}

\begin{proof}[Proof of Lemma~\ref{l:Rfarcyl}]$^\dagger$
{\it Step 1: canonical neighborhoods.} On the graphical part of the end, $R_0\sqrt\tau\le\rho\le\hat R(s)\sqrt\tau$, Lemma~\ref{l:Rss}(i) and the interior estimates of \cite{EH2} make the flow close at the scale $\sqrt\tau$ to the cylindrical end of $S$, where $H_S\to1/\sqrt2$. So, once $R_0$ is large and $r_1$ small, it is mean convex there with $H\sqrt\tau\asymp1$, and $\varepsilon_2$-close to a cylinder at  scale $1/H$.

For the far part, let $y_0$ be the point of the axis at distance $r\le4r_1$ from $x_0$, and let $s_*=s_{2r}$ and $\sqrt{\tau_*}=2r/\hat R(s_*)$. So $y_0$ is at rescaled distance $\hat R(s_*)/2$ at the time $s_*$, and in the graphical region before it. By Lemma~\ref{l:Rss}(i),(iii), the asymptotics of the end \cite[Theorem~1.10]{BK} and Corollary~\ref{c:BK1}, at the times $s\in[s_*-2\log\hat R(s_*),s_*]$ the flow near $y_0$ is, at the scale $\sqrt{\tau_*}$, $\delta_*$-close to the round cylinder that pinches at time $T$, on a region of rescaled size comparable to $\hat R(s_*)$, cf.\ \cite{CM1}. Here $\delta_*=C\varepsilon(s_*)+\eta_c(\hat R(s_*))$, where $\eta_c(\Xi)\to0$ measures how close the end of $S$ is to its asymptotic cylinder at distance $\Xi$, so $\delta_*\to0$ as $r\to0$. Thus, for $r_1$ small, Lemma~\ref{l:RCHHq}, with $L_Q=2C$ and its $\epsilon$ small depending on $\varepsilon_2$, applies at $(y_0,t_*)$. It gives $H>0$ and canonical neighborhoods (a sphere, a cylinder, a bowl or an ancient oval) on $B_{2C\sqrt{\tau_*}}(y_0)\times[t_*,T)$. A point $(x,t)$ of $F$ with $|x-x_0|=r\ge\hat R(s)\sqrt\tau$ lies in this region: $t\ge t_*$, and $\tau\le r^2/\hat R(s)^2\le\tau_*/4$ since $\hat R$ is nondecreasing, so $x$ is within $\frac C2\sqrt{\tau_*}$ of the axis, hence within $C\sqrt{\tau_*}$ of $y_0$. The other points of $F$ lie on the graphical part. Finally, an ancient oval of large rescaled diameter is close on the ball to a bowl or a cylinder, since such ovals, rescaled at points with $|A|=1$, subconverge to a noncompact, non-flat ancient low entropy flow \cite[Theorem~1.2]{ChHH}.

{\it Step 2: $H\geq cR_0/\rho$.} On the graphical part this follows from $H\sqrt\tau\asymp1$ and $\rho\ge R_0\sqrt\tau$. On the far part, where $t\ge t_{2\rho}$, Lemma~\ref{l:Rsheets} puts $\Sigma_t$ near the end within $C\rho/\hat R(s_{2\rho})\le C\rho/R_0$ of the axis, for $r_1$ small. Every model of Step~1 has width at least $c/H$ in every direction at the scale $1/H$, and the canonical neighborhood must fit there, so $H\ge cR_0/\rho$.

{\it Step 3: no compact models.} If the canonical neighborhood at a point $p$ of $\Sigma_t\cap F$ were a sphere or an ancient oval of rescaled diameter at most $D$, the component of $\Sigma_t$ through $p$ would be a closed surface of diameter at most $CD\rho/(cR_0)$, much smaller than $\rho$, and so not the component of the graph in Lemma~\ref{l:Rss}(i). Since $|p-x_0|\le4r_1$ and $T-t\le(4r_1/R_0)^2$, this contradicts Lemma~\ref{l:Rss}(v) once $r_1$ is small. This proves (1).

{\it Step 4: (2).} On $F$, $\rho=|x-x_0|$, so by Lemma~\ref{new:radial}(1), for a function $f(\rho,\tau)$,
\begin{align}
	(\partial_t-\cL_{H^2})f=-\partial_\tau f+\partial_\rho f\,\frac{|\nabla^T\rho|^2-2}\rho-\partial_\rho^2f\,|\nabla^T\rho|^2-2\,\partial_\rho f\,\big\langle\nabla\log H,\nabla^T\rho\big\rangle\, ;\notag
\end{align}
the last term is the drift. Write $\Psi_\pm=\rho^{1-2\alpha\mp2\theta}$, so that $B_\pm=P_\pm\Psi_\pm$ with $P_\pm=1\pm\eta\mp\eta cR_0\,\rho^{-1}H^{-1}$. By (1), $\rho H\geq cR_0$, which gives the bounds for $P_\pm$ in (2), and in particular $\frac12\le P_\pm\le\frac32$.

Since $H$ is a Jacobi field, the quotient rule of Section~\ref{ss:Rjacobi} gives $(\partial_t-\cL_{H^2})H^{-1}=-|A|^2/H\leq-H/2$. 
This term has the right sign in both cases and it will dominate the other terms; we will use here that $H$ is positive and decays at worst like $\rho^{-1}$ by (1). 
By the formula above, $|(\partial_t-\cL_{H^2})\rho^{-1}|\leq4\rho^{-3}+2C_AH\rho^{-2}$. Also $|\nabla\rho^{-1}|\leq\rho^{-2}$ and $|\nabla H^{-1}|=|\nabla H|/H^2\leq C_A$. By the product rule,
\begin{align}
	\pm(\partial_t-\cL_{H^2})P_\pm=-\eta cR_0\,(\partial_t-\cL_{H^2})(\rho^{-1}H^{-1})\geq\eta cR_0\Big(\frac{|A|^2}{\rho H}-\frac4{\rho^3H}-\frac{4C_A}{\rho^2}\Big)\geq\frac{\eta cR_0}4\,\frac H\rho\, .\notag
\end{align}
The last step uses $|A|^2\geq H^2/2$ and $\rho H\geq cR_0$, and holds once $cR_0\geq32C_A+6$.

 Since $0<\theta\le\frac14$, $|\partial_\rho\log\Psi_\pm|\le2\alpha/\rho$ and $\rho^2|\partial_\rho^2\Psi_\pm|\le C_2\Psi_\pm$, with $C_2$ depending only on $\alpha$. So, by the formula above,
\[
|(\partial_t-\cL_{H^2})\Psi_\pm|\le\Psi_\pm\Big(\frac{C_2}{\rho^2}+\frac{4\alpha C_AH}\rho\Big).
\]
Also $|\nabla P_\pm|\le\eta cR_0(\rho^{-2}H^{-1}+\rho^{-1}C_A)\le\eta(1+cR_0C_A)/\rho$ and $|\nabla\Psi_\pm|\le2\alpha\Psi_\pm/\rho$. By the product rule, using $\frac12\le P_\pm\le\frac32$ and $\rho^{-2}\le H/(cR_0\rho)$,
\[
\pm(\partial_t-\cL_{H^2})B_\pm\ \ge\ B_\pm\,\frac H\rho\Big(\frac{\eta cR_0}6-C_2-8\alpha C_A-4\alpha\Big).
\]
Let $K_P=\max\{6(C_2+8\alpha C_A+4\alpha+1+4C_A),\,32C_A+6\}/c$. If $\eta R_0\ge K_P$, the bracket is at least $1+4C_A$, which gives the inequality of (2); and since $\eta\le\frac12$, $R_0\ge2K_P$, so $cR_0\ge cK_P\geq32C_A+6$, as used above. Finally, $(\partial_t-\Delta-|A|^2)(fH)=H(\partial_t-\cL_{H^2})f$ by the quotient rule of Section~\ref{ss:Rjacobi}, so $HB_+$ is a supersolution and $HB_-$ a subsolution of the Jacobi equation on $F$.
\end{proof}

\section{The bowl soliton}	\label{app:bowl}\label{ss:Rbowlineq}

 Near the tips of degenerate necks,
  the flow is close to a rotationally symmetric translating solution known as the bowl soliton. 
 We use the following facts about the bowl in Lemma~\ref{l:Rfarcyl} and in Section~\ref{ss:Rbarrier}. The unit-speed bowl in $\RR^3$ is the graph of $f(r)$ over a plane perpendicular to the direction $e$ of translation.  The translator equation
 $H= -\langle e , \nn \rangle$ gives that, 
 after subtracting a constant, 
  $f$  satisfies (cf. section $6$ in \cite{BCh})
\begin{align}
	\frac{f''}{1+f'^2}+\frac{f'}r=1\, , f(0) = f' (0) = 0 {\text{ and }} f''(0) = \frac{1}{2}  .\notag
\end{align}
Taylor expansion gives near zero that $f(r) = \frac{r^2}{4} + O(r^4)$ and $f' (r) = \frac{r}{2} + \frac{r^3}{32} + O(r^4)$.

If we define $p$ and $D$ by  $ p=\frac{f'}{r}$ and $ D=r-f'$, then 
  $H=(1+(f')^2)^{-1/2}$, the principal curvatures are $(1-p)H$ and $pH$, and $|\nabla\log H|=f'(1-p)H=pDH$. 
  The ratio of $|\nabla \log H|$ and $|A|$ determines the sign of the zero-th order term in the Jacobi-type equation
  for powers of $H$.  To keep track of this, define
\begin{align}	\label{e:RRbowl}
	R:=\frac{|\nabla\log H|^2}{|A|^2}=\frac{p^2D^2}{(1-p)^2+p^2}\leq D^2\, .
\end{align}

\begin{Lem}	\label{l:Rbowlineq}
On the bowl:
\begin{enumerate}
\item $p$ is strictly increasing, with $\frac12<p<1$. In particular, $f' , f'' > 0$ for $r>0$.
\item $D<1$, so $R<1$; moreover $D\to0$ as $r\to\infty$, so $R_{\max}:=\sup R<1$.
\end{enumerate}
\end{Lem}

\begin{proof}
(1) From the equation, $rp'=\Psi(r,p):=(1+r^2p^2)(1-p)-p$. Since $\Psi(r,\frac12)=\frac{r^2}8>0$ and $\Psi(r,1)=-1<0$, the function $p$, which starts at $\frac12$, stays in $(\frac12,1)$. Near $r=0$, $p=\frac12+\frac{r^2}{32}+O(r^4)$, so $p'>0$. At a first zero $r_*$ of $p'$ we would have $p''(r_*)\leq0$, but differentiating $rp'=\Psi$ gives $r_*p''(r_*)=\partial_r\Psi=2r_*p^2(1-p)>0$.

(2) For $r\leq2$, $f'=rp>r-1$, and thus $D<1$ on $[0,2]$ (note that $D(0)=0$). To see that $D<1$ for all $r>2$, suppose instead that there is some first $r'>2$ with $D(r') = 1$. It follows that $D'(r') \geq 0$, so that
$$
	0 \leq 1 - f''(r') = 1 - \frac{D(r')}{r'} \left( 1 + ( f'(r' ))^2 \right)= \frac{ - r'^2 + 3\, r' - 2}{r'}
	= -\frac{ (r' -1)(r' - 2)}{r'}
\, .
$$
This is impossible for $r_0 >2$, so we conclude that $D<1$ as claimed. We show next that $D \to 0$ as $r\to \infty$.
Note that 
 $D'=1-(1+(r-D)^2)D/r\leq1-(r-2+\frac2r)D$ for $r\geq1$. Given $c_0\in(0,1)$, for $r>2+2/c_0$ the right side is at most $-1$ wherever $D\geq c_0$, so $D$ eventually stays below $c_0$. Since $R$ is continuous, $R<1$, and $R\to0$ at both ends, its supremum is attained and is less than $1$.
\end{proof}

\end{document}